\documentclass[12pt]{amsart}

\usepackage{amsmath, amsthm, amssymb}
\usepackage{graphicx}
\usepackage{mathbbol}
\usepackage{txfonts}
\usepackage[usenames]{color}
\usepackage{pgfplots}
\usepackage{comment}

\allowdisplaybreaks[4]

\newtheorem{thm}{Theorem}[section]
\newcommand{\bt}{\begin{thm}}
\newcommand{\et}{\end{thm}}

\newtheorem{conj}[thm]{Conjecture}

\newtheorem{cor}[thm]{Corollary}   
\newcommand{\bc}{\begin{cor}}
\newcommand{\ec}{\end{cor}}

\newtheorem{lem}[thm]{Lemma}   
\newcommand{\bl}{\begin{lem}}
\newcommand{\el}{\end{lem}}

\newtheorem{prop}[thm]{Proposition}
\newcommand{\bp}{\begin{prop}}
\newcommand{\ep}{\end{prop}}

\newtheorem{defn}[thm]{Definition}
\newcommand{\bd}{\begin{defn}}    
\newcommand{\ed}{\end{defn}}

\newtheorem{rmrk}[thm]{Remark}   

\newcommand{\br}{\begin{rmrk}}
\newcommand{\er}{\end{rmrk}}

\newtheorem{example}[thm]{Example}

\newcommand{\MinA}{\operatorname{MinA}}

\newcommand{\Scal}{\operatorname{Scalar}}

\newcommand{\be}{\begin{equation}}

 \newcommand{\ee}{\end{equation}}

\newcommand{\R}{\mathbb{R}}

\newcommand{\diam}{\operatorname{Diam}}

\newcommand{\area}{\operatorname{Area}}
\newcommand{\vol}{\operatorname{Vol}}

\newcommand{\Sph}{{\mathbb S}}         

\newcommand{\og}{\overline{\Gamma}}
\newcommand{\on}{\overline{\nabla}}

\begin{document}

\title[An Extreme Limit with Nonnegative Scalar Curvature]{An Extreme Limit with Nonnegative Scalar}

\author{Christina Sormani}
\thanks{Prof. Sormani was partially supported by NSFDMS1006059 and PSC CUNY Funding.}
\address{CUNY Graduate Center and Lehman College, NY, NY, 10016}
\email{sormanic@gmail.com}

\author{Wenchuan Tian}
\address{Department of Mathematics, University of California, Santa Barbara, CA 93106-3080}
\email{tian.wenchuan@gmail.com}

\author{Changliang Wang}
\thanks{Changliang Wang was partially supported by the Fundamental Research Funds for the Central Universities.} 
\address{School of Mathematical Sciences and Institute for Advanced Study, Tongji University, Shanghai 200092, China}
\email{wangchl@tongji.edu.cn}


\begin{abstract}
In 2014, Gromov vaguely conjectured that a sequence of manifolds with nonnegative scalar curvature should have a subsequence which converges in some weak sense to a limit space with some generalized notion of nonnegative scalar curvature.  The conjecture has been made precise at an IAS Emerging Topics meeting: requiring that the sequence be three dimensional with uniform upper bounds on diameter and volume, and a positive uniform lower bound on MinA, which is the minimum area of a closed minimal surface in the manifold.  Here we present a sequence of warped product manifolds with warped circles over standard spheres, that have circular fibres over the poles whose length diverges to infinity, that satisfy the hypotheses of this IAS conjecture.  We prove this sequence converges in the $W^{1,p}$ sense for $p<2$ to an extreme limit space that has nonnegative scalar curvature in the distributional sense as defined by Lee-LeFloch and that the total distributional scalar curvature converges.   This paper only requires expertise in smooth Riemannian Geometry, smooth minimal surfaces, and Sobolev Spaces.  In a second paper, requiring expertise in metric geometry, the first two authors prove intrinsic flat and Gromov-Hausdorff convergence of our sequence to this extreme limit space and investigate its geometric properties. 
\end{abstract}

\maketitle

\section{Introduction}

For nearly half a century, mathematicians have been studying classes of Riemannian manifolds with various curvature bounds, whether sequences of spaces in these classes have converging subsequences, and what kinds of properties hold on their limit spaces.   With strong enough curvature and injectivity radius bounds one can obtain smooth convergence of the Riemannian manifolds, and with less strong hypotheses one obtains $C^{1,\alpha}$ convergence to a smooth manifold with a $C^{1,\alpha}$ metric tensor.  With only nonnegative sectional curvature, one has Gromov-Hausdorff convergence to an Alexandrov space with nonnegative curvature, and with only Ricci curvature, one has Gromov-Hausdorff Convergence to a connected geodesic metric space called an RCD space.    See the texts of Petersen \cite{Petersen-text} and Gromov \cite{Gromov-metric} and the survey by the first author \cite{Sormani-Scalar}.

In \cite{Gromov-Plateau} and \cite{Gromov-Dirac}, Gromov conjectured that a sequence of Riemannian manifolds with nonnegative scalar curvature, $\Scal \ge 0$, should have a subsequence which converges in some weak sense to a limit space with some generalized notion of ``nonnegative scalar curvature".   This conjecture was made more precise at an IAS Emerging Topics Workshop co-organized by Gromov and the first author as follows \cite{Sormani-Scalar}:

\begin{conj}\label{Scalar-Compactness}
Let $\{M_j^3\}$ be a sequence of
 closed oriented three dimensional manifolds without boundary satisfying
\be \label{other-three}
\Scal_j \ge 0,
\qquad 
 \vol(M_j) \le V,
\qquad
 \diam(M_j) \le D,
\ee
and
\be \label{defn-MinA}
\MinA(M^3_j)=\inf\{\area(\Sigma)\, : \, \Sigma \textrm{ closed min surf in } M_j^3\,\} \ge A_0>0.
\ee
Then a subsequence of $\{M_j\}$ converges in the volume preserving intrinsic flat sense to a three dimensional rectifiable limit space, $M_\infty$.
Furthermore, $M_\infty$, is a connected geodesic metric space, that has Euclidean tangent cones almost everywhere, and has nonnegative generalized scalar curvature.   
\end{conj}

The $\MinA$ condition in (\ref{defn-MinA}) can be viewed as a noncollapsing condition which prevents counter examples like sequences of round spheres rescaled to a point.  
Note that in a warped spheres with 
\be \label{warped-spheres}
g_j=dr^2+h_j(r)^2 g_{\Sph^2} \textrm{ where } h_j:[0,D_j]\to (0,\infty)
\ee
any level set $r^{-1}(r_0)$ where $h_j'(r_0)=0$ is a minimal surface of area
$4\pi h_j(r_0)^2$.   In joint work with Jiewon Park, the second and third authors have proven Conjecture~\ref{Scalar-Compactness}
for sequences of warped spheres with metric tensors of the form (\ref{warped-spheres})
\cite{Park-Tian-Wang-18}.   In particular they show the $\MinA$ condition prevents the formation of a thin tunnel between two noncollapsed regions.  More precisely, they
prove a subsequence of the $h_j$ converge
in the $C^0$ and $H^1_{loc}$ sense to a bounded function
$h_\infty:[0,D_\infty]\to [0,\infty)$.  Although $h_\infty$ can equal zero, the set 
$\{r:\, h_\infty(r)>0\}$ is a connected set.  

The $\MinA$ condition in (\ref{defn-MinA}) was added to the conjecture in light of sequences of $M_j^3$ satisfying the other three hypotheses (\ref{other-three}) of this conjecture constructed by Basilio, Dodziuk, and the first author whose limit spaces do not satisfy the
properties of spaces with a natural notion of generalized nonnegative scalar curvature \cite{BDS-sewing}.  These sequences have increasingly many increasingly tiny tunnels.   One of their limit spaces is a standard round sphere with a great circle collapsed to a point.    More recently Basilio and the first author constructed sequences
of $M_j^3$ satisfying (\ref{other-three}) whose limit space is a three-sphere with an arbitrary region collapsed to a point \cite{BS-seq}.  Basilio, Kazaras, and Sormani constructed a sequence of $M_j^3$ satisfying  (\ref{other-three}) whose limit space is a metric space with no geodesics \cite{BKS-no-geod}.  In higher dimensions even more
problematic issues can occur as seen in the work of Lee-Naber-Neumeyer \cite{Lee-Naber-Neumayer-dp} and Lee-Topping \cite{Lee-Topping-metric}.

In Section~\ref{Sect-Seq}, we present a new example: a sequence of
$M_j^3$ diffeomorphic to $\Sph^2\times \Sph^1$, with metric tensors of the
form
\be\label{g-warp}
g_j=g_{\Sph^2}+f_j^2 g_{\Sph^1} \textrm{ where } f_j:\Sph^2\to (0,\infty)
\ee
that satisfy all the hypotheses of Conjecture~\ref{Scalar-Compactness} [Example~\ref{ex-sequence}].
This new example is of interest because the $f_j$ diverge to infinity above the poles, stretching two circular fibres to infinite length as $j\to \infty$ [Lemma~\ref{lem-fibres}].
Despite this stretching in one direction, we prove the diameters 
and the volumes are uniformly bounded in Proposition~\ref{UniformBoundDiameter}
and Proposition~\ref{prop-vol-sequence} respectively.  

In Section~\ref{Sect-conv}, we introduce an extreme warped product space, $\Sph^2\times \Sph^1$ with
a metric tensor $g_\infty$ of the form described in (\ref{g-warp}) with a warping factor 
$f_\infty: \Sph^2\to (0,\infty]$ that is equal to $\infty$ at the poles 
so that $g_\infty$ is a smooth metric tensor away from a singular set $S$ which consists of the two circular fibres above the poles [Example~\ref{ex-limit}].   In Proposition~\ref{prop-warp-conv} we prove our sequence of $g_j$ of Example~\ref{ex-sequence} converge to $g_\infty$ smoothly away $S$.  In Proposition ~\ref{metric-W1p} we prove
$g_\infty$ is in $W^{1,p}(\Sph^2\times \Sph^1)$ for $p\in [1,2)$ 
but not $H^1_{loc}$.   

We postpone all discussion of the distance functions on this extreme limit space and how to view this extreme limit space as a metric space to a separate paper
by the first two authors \cite{ST-geom-ex}.   That paper also contains a proof that our 
sequence of Riemannian manifolds $M_j$ converges in the intrinsic flat and Gromov-Hausdorff sense to our extreme warped product limit space.  Thus it is essential that any
notion of generalized nonnegative scalar curvature being considered for the statement of Conjecture~\ref{Scalar-Compactness} must hold on our limit space.

In Section~\ref{Sect-curv} we study the curvature of our extreme warped product space, $(\Sph^2\times \Sph^1, g_\infty)$, because it can be used to test notions of generalized nonnegative scalar curvature.  It is a particularly interesting example, because the singularities are not isolated.   Note that smooth metric tensors with nonnegative scalar curvature and isolated singularities have been studied by 
Korevaar, Mazzeo, Pacard, and Schoen in \cite{Kor-Maz-Pac-Schoen}.  More
recently, Li and Mantoulidis have studied manifolds with nonnegative curvature and skeleton singularities in \cite{Li-Mantoulidis-skeleton}, however they do not allow the
metric tensor to blow up along the edge singularities as ours does.

We calculate the Ricci curvature of $g_\infty$ where it is smooth and see that it has no lower bound in Lemma~\ref{LemmaRicci}.  As a consequence our space cannot be a CD space as defined
in work of Lott-Villani \cite{Lott-Villani} and Sturm \cite{Sturm}, nor an RCD space as defined in the work of 
Ambrosio-Gigli-Savare \cite{Ambrosio-Gigli-Savare}.   See Remark~\ref{rmrk-not-RCD}.
It is possible that our calculation might be useful for those who would like to study the Ricci flow emanating from $g_\infty$.  Recall that Bamler used Ricci flow to prove that when a sequence of smooth manifolds with nonnegative scalar curvature converge in the $C^0$ sense to a smooth manifold, then the limit manifold has nonnegative scalar curvature as well \cite{Bamler-C}.  Burkhardt-Guim \cite{B-G-GAFA} then used Ricci flow to define a generalized notion of nonnegative scalar curvature for $C^0$ metric tensors in \cite{B-G-GAFA}.   See Remark~\ref{rmrk-Ricci-flow} for a discussion how one might extend this to our even lower regularity metric tensor.

In Theorem~\ref{thm-distr-scalar}, we prove the limit space, $(\Sph^2\times \Sph^1, g_\infty)$ has nonnegative distributional scalar curvature in the sense of Lee-LeFloch  
\cite{Lee-LeFloch}.   Lee-LeFloch defined their notion building upon the work of LeFloch-Mardare \cite{LM07} and they were able to prove the positive mass theorem for manifolds whose metric tensors are $C^0\cap W^{1,n}$ using this notion. See also the
work of Cecchini, Simone and Hanke \cite{CHS-Lip} and the work of
Lee and Tam \cite{Lee-Tam-rigidity} for other application of Lee-LeFloch's notion of distributional scalar curvature.   Although our metric tensors are not as regular as the
metric tensors they studied, we are nevertheless able to prove we have nonnegative
scalar curvature in their sense [Remarks~\ref{rmrk-Lee-LeFloch-original}-\ref{rmrk-Lee-LeFloch}].   One may view our example as yet another validation of the Lee-LeFloch notion.   

In Definition \ref{defn-distr-total-scalar-curvature} we extend the notion of total distributional scalar curvature to a natural one defined in the sense of Lee-LeFloch.   In Lemma~\ref{lem-distr-total-scalar-curvature}, we prove the total scalar curvature of 
our sequence, $M_j=(\Sph^2\times \Sph^1,g_j)$, 
converges to the total distributional scalar curvature
of our extreme warped limit space, $M_\infty=(\Sph^2\times \Sph^1,g_\infty)$.   In
Remark~\ref{rmrk-total-scalar} we observe that although the
singular set consists of only two fibres it contributes a positive amount to this total distributional scalar curvature.   See the work of Hamanaka exploring the behaviour of the total scalar curvature under $C^0$ convergence in \cite{Hamanaka-total}.

Finally in Section~\ref{sect-MinA}, we prove our sequence satisfies the $\MinA$ hypothesis of Conjecture~\ref{Scalar-Compactness}.
Note that it is incredibly difficult to find a uniform lower bound on $\MinA$ for a Riemannian manifold.  Although we believe that the manifolds $M^3_j$ our sequence in Example~\ref{ex-sequence}, have $\MinA(M^3_j)\ge 4\pi$ (see Conjecture~\ref{Conj-MinA}), we were not able to prove this.  Instead, in Theorem~\ref{thm-MinA}, we complete a proof by contradiction that there is a uniform lower bound on $\MinA(M_j)$ for our sequence using the smooth convergence away from the singular set and the mean curvature of the boundaries of regions which have smooth convergence.  This work and additional lemmas exploring the areas of key minimal surfaces in our sequence in appears in Section~\ref{sect-MinA}.
We believe that the techniques we apply there should be useful in other settings as well.  

\vspace{1cm}
\noindent
{\bf Acknowledgements:}

We would like to thank the IAS for hosting the Emerging Topics Working Group on Scalar Curvature and Convergence in 2018. Christina Sormani was partial funded by NSF DMS 1612409 and a PSC CUNY award. Changliang Wang was partially supported by the Fundamental Research Funds for the Central Universities.  We are all grateful to the Fields Institute for hosting the {\em Summer School on Geometric Analysis} in July 2017 where Christina Sormani gave a series of lectures and met Wenchuan Tian and Changliang Wang.

\section{The Sequence of Riemannian Manifolds}
\label{Sect-Seq}

In this section we define a sequence of Riemannian manifolds and prove they
satisfy the uniform bounds in the hypotheses of Conjecture \ref{Scalar-Compactness}.

\begin{figure}[h] 
   \centering
   \includegraphics[width=4in]{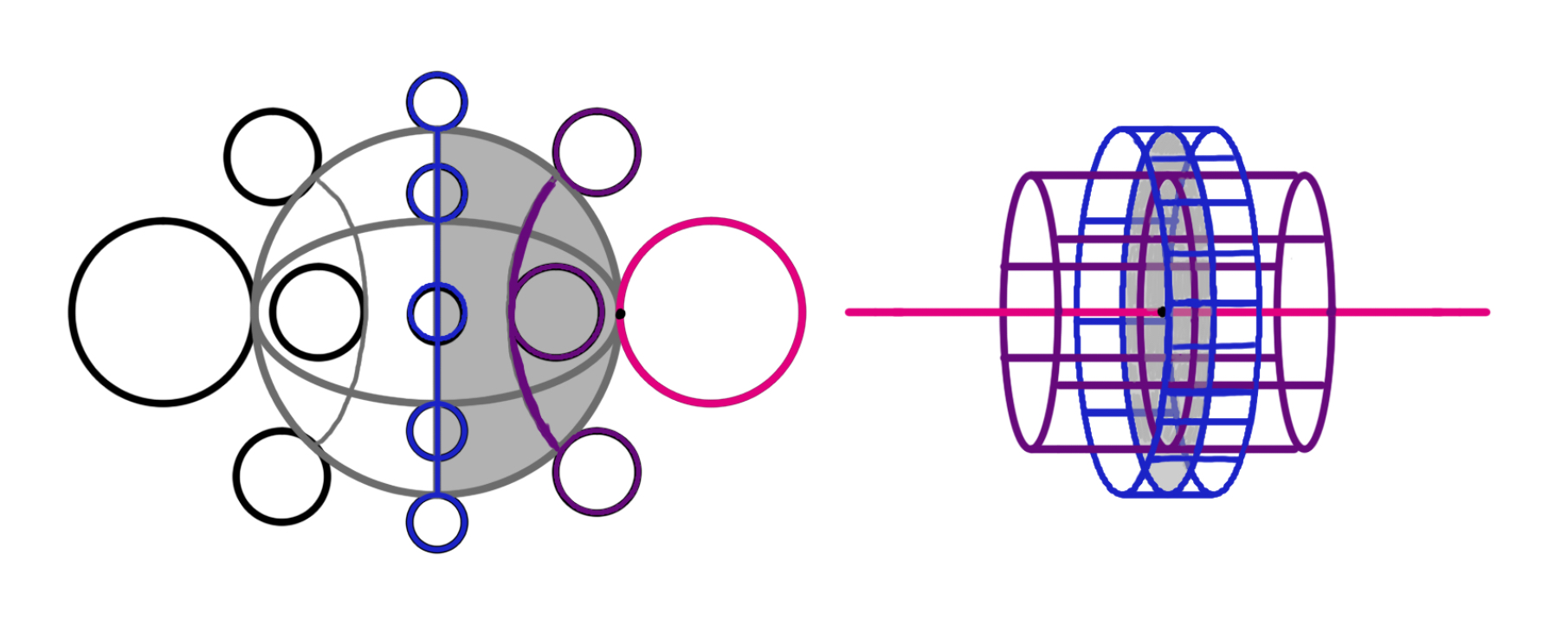} 
   \caption{On the left we have $\Sph^2\times_{f_j}\Sph^1$ with circles above every point in $\Sph^2$
   of varying size depending on $f_j$.   On the right, we look at one hemisphere of $\Sph^2$
   viewed as a disk, with the circles of varying size depicted as horizontal line segments of length $2\pi f_j(r)$.  
   }
   \label{fig-STW-1}
\end{figure}

\begin{example}\label{ex-sequence}
	Consider the sequence of warped product Riemannian manifolds, $\Sph^2\times_{g_j} \Sph^1$, which are diffeomorphic to $\Sph^2\times \Sph^1$, with Riemannian metrics 
	\be
	g_j=g_{\Sph^2}+f_j^2(r) g_{\Sph^1}=dr^2+ \sin^2(r) d\theta^2+f_j^2(r) d\varphi^2.
	\ee
	written using $(r,\theta)$ coordinates on $\Sph^2$ and $\varphi$ on the $\Sph^1$-fibres
	where 
	\be
	f_j(r)=\ln\left(\frac{1+a_j}{\sin^2 r+a_j}\right)+\beta,
	\ee
taking $\beta\geq 2$, positive $a_j$ decreasing to $0$.   See Figure~\ref{fig-STW-1}.
 \end{example}
 
In Subsection~\ref{subsect-smooth} we confirm that $\Sph^2\times_{g_j} \Sph^1$
defined in Example~\ref{ex-sequence}  are indeed smooth Riemannian
manifolds with no singularity at $r=0$ and $r=\pi$.   
We show the functions, $f_j$, viewed as functions on $\Sph^2\times_\Sph^1$ are smooth
in Lemma~\ref{lem-sin2}, increasing with $j$ in Lemma~\ref{lem-warp-inc}, and unbounded as $j\to \infty$ in Lemma~\ref{lem-warp-unbounded}.

In Subsection~\ref{subsect-scalar} we prove $\Sph^2\times_{g_j} \Sph^1$
defined in Example~\ref{ex-sequence} have nonnegative scalar curvature [Proposition~\ref{prop-scalar-seq}].  

In Subsection~\ref{subsect-vol} we prove a uniform upper volume bound on the sequence 
of Riemannian manifolds in Proposition~\ref{prop-vol-sequence} and also control the areas of the level sets $r^{-1}(r_0)$ in these spaces.

In Subsection~\ref{subsect-fibres} we prove the fibres above the poles have lengths diverging to infinity in Lemma~\ref{lem-fibres}.

In Subsection~\ref{subsect-unif-dist} we view the Riemannian manifolds $(\Sph^2\times\Sph^1,g_j)$ as metric spaces with distance functions, $(\Sph^2\times\Sph^1,d_j)$.  We prove 
the distance functions are increasing in Lemma~\ref{LemDistanceFunctionsInc} and prove
uniform bounds on the distance functions in Lemma~\ref{LemmaBoundDistance}.   We also prove the
diameter is uniformly bounded from above (even though the metric tensors are not bounded) in Proposition~\ref{UniformBoundDiameter}.

\subsection{The Sequence is Smooth but Increasing and Unbounded}
\label{subsect-smooth}

Here we prove functions and metric tensors in Example \ref{ex-sequence} are smooth in Lemma~\ref{lem-sin2}.
We show they are increasing in Lemma~\ref{lem-warp-inc} and unbounded in Lemma~\ref{lem-warp-unbounded}.

\begin{lem}[Regularity of the Metric] \label{lem-sin2} 
The function $f_j$ and the metric tensor $g_j$ are
smooth in Example \ref{ex-sequence} .
\end{lem}

\begin{proof}
The function $\sin^2(r)$ is smooth on $\Sph^2$
and thus $f_j(r)$ is smooth on $\Sph^2$.  Thus the metric $g_j$ is smooth
as a metric tensor on $\Sph^2\times\Sph^1$.   For more details see below.

Denote $N$ as the north pole in $\Sph^2$, we use the coordinate $(r,\theta)$ to denote the polar coordinate centered at $N$. In this coordinate the metric in $\Sph^2\times_{f_j}\Sph^1$ is $g_j=dr^2+\sin ^2 r d\theta^2+f_j^2(r)d\varphi^2$. 
	
Consider the point $p\in \Sph^2$ where $r(p)=\frac{\pi}{2}$ and $\theta(p)=0$. We use $(\tilde{r},\tilde{\theta})$ to denote the polar coordinate centered at the point $p$, then in the new coordinate the metric $g_j=d\tilde{r}^2+\sin^2 \tilde{r}d\tilde{\theta}^2+\tilde{f}_j^2(\tilde{r},\tilde{\theta})d\tilde{\varphi}^2$, where we choose the geodesic in $\Sph^2$ connecting $N$ and $p$ as the line for $\tilde{\theta}=0$.

	 Consider a point $q\in \Sph^2$ such that $q=(r,\theta)$ in the polar coordinate centered at $N$ and that $q=(\tilde{r},\tilde{\theta})$ in the polar coordinate centered at $p$. By the cosine law in the sphere, we have
	 \be
	 \cos r= \cos \frac{\pi}{2}\cos\tilde{r}+\sin \frac{\pi}{2} \sin \tilde{r} \cos\tilde{\theta},
	 \ee
	 and as a result,
	 \be
	 \sin^2 r=1-\sin^2\tilde{r}\cos^2\tilde{\theta}.
	 \ee
	 Hence
	 \be
	 \widetilde{f}_j(\tilde{r},\tilde{\theta})=\ln\left(\frac{1+a_j}{1-\sin^2\tilde{r}\cos^2\tilde{\theta}+a_j}\right)+\beta.
	  \ee
	  As a result we conclude that $f_j$ as a function on ${\mathbb S}^2$ is smooth
	  and so the metric tensor, $g_j$, is smooth on $\Sph^2\times\Sph^1$.
\end{proof}

\begin{lem} \label{lem-warp-inc}
For $f_j(r)$ and $g_j$ as defined in Example \ref{ex-sequence} we have
\be
f_j(r) \le f_{j+1}(r) \textrm{ and } g_j \le  g_{j+1}.
\ee
\end{lem}

\begin{proof}
This follows because $a_j$ are decreasing and $\ln$ is increasing and
\be
\frac{d}{da} \,\left(\frac{1+a}{\sin^2 r+a} \right)
=\frac{ (1)(\sin^2r+a) - (1+a)(1)}{(\sin^2 r+a)^2}\le 0.
\ee
\end{proof}

\begin{lem} \label{lem-warp-unbounded}
The sequence $f_j$ as defined in Example \ref{ex-sequence} is unbounded in the sense that
\be
\lim_{j\to \infty} f_j(0) \to \infty \textrm{ and } \lim_{j\to \infty} f_j(\pi) \to \infty.
\ee
\end{lem}

\begin{proof}
We have 
\be
f_j(0)=f_j(\pi)=\ln\left(\frac{1+a_j}{0+a_j}\right)+\beta
\ee
and $a_j$ decrease to $0$, so $(1+a_j)/a_j \to \infty$.
\end{proof}

\subsection{Proving the Sequence has Nonnegative Scalar Curvature:}
\label{subsect-scalar}

Here we prove the sequence in  Example~\ref{ex-sequence} has nonnegative scalar curvature.

As the curvatures of warped product manifolds have been computed in the past, we simply consulted the article of Dong-Soo Kim and Young Ho Kim \cite{KK-compact} rather than completing the computation from scratch.   In Proposition 2 of \cite{KK-compact}, the scalar
curvature of a warped product metric tensor of the form 
\be
g=g_{{\mathbb S}^2}+ f^2(u)d\varphi^2 \textrm{ where } f:{\mathbb S}^2 \to {\mathbb R}^+
\ee
was computed to be
\be
\Scal=2-2\frac{\Delta f}{f} =2 \left(\frac{f - \Delta f }{f}\right) 
\ee
where  $\Delta$ is the Laplace Beltrami operator on the standard round sphere, 
$({\mathbb S}^2, g_{{\mathbb S}^2})$.   This
particular partial differential equation has been studied in relation to scalar curvature and minimal surfaces in the past.  See, for example, the work of Schoen-Yau \cite{Schoen-Yau-structure} \cite{Schoen-Yau-singularities} and Kai Xu \cite{xu-geometric}.  

Note that in particular
\be \label{eq-scal-equiv}
\Scal \ge 0 \iff f \geq \Delta f.
\ee
In fact, we studied this equation for some time while trying to prove 
Conjecture~\ref{Scalar-Compactness}.
in this setting before devising our example.

\begin{prop}[Scalar Curvature] \label{prop-scalar-seq}
	The manifold $\Sph^2\times_{f_j}\Sph^1$ as defined in Example \ref{ex-sequence} has 
	\be
\Scal_j = 2-2\frac{\Delta f_j}{f_j} =  2 \left(\frac{f_j - \Delta f_j }{f_j}\right) \ge 0.
\ee
\end{prop}

\begin{proof}
In polar coordinates in $\Sph^2$ we have
\be
\Delta f_j=\partial_r^2 f_j+\frac{\cos r}{\sin r} \partial_r f_j 
\ee
because $f_j$ only depends on $r$.
We have
\be
 \partial_r f_j= -\frac{2 \cos r \sin r}{\sin^2 r+a_j},
\ee
and
\be
 \partial_r^2 f_j=\frac{4 \cos^2 r \sin^2 r}{(\sin^2 r+a_j)^2} -\frac{2\cos^2 r}{\sin^2 r+a_j}+\frac{2\sin^2 r}{\sin^2 r+a_j}.
\ee
As a result 
\be
\begin{split}
\Delta f_j &=\frac{4 \cos^2 r \sin^2 r}{(\sin^2 r+a_j)^2} -\frac{4\cos^2 r}{\sin^2 r+a_j}+\frac{2\sin^2 r}{\sin^2 r+a_j}\\
&=\frac{-4a_j \cos^2 r + 2a_j\sin^2 r +2\sin^4 r}{(\sin^2 r+a_j)^2} \\
&=\frac{2(\sin^2r +a_j)^2 -2a_j\sin^2 r-2a_j^2-4a_j \cos^2 r }{(\sin^2 r+a_j)^2} \\
&\leq 2.
\end{split}
\ee
On the other hand, since $\beta\geq 2$ and $1+a_j\geq \sin^2 r+a_j$ we have
\be
f_j(r)=\ln\left(\frac{1+a_j}{\sin^2 r+a_j}\right)+\beta\geq 2.
\ee
Thus $f_j\geq \Delta f_j$ and we are done.
\end{proof}

\subsection{Proving a Uniform Volume Bound for the Sequence:}
\label{subsect-vol}

Here we prove a uniform volume bound on the sequence 
of Riemannian manifolds defined in Example \ref{ex-sequence} in
Proposition~\ref{prop-vol-sequence}. We
 also control the areas of the level sets $r^{-1}(r_0)$ in these spaces
 in Lemma~\ref{lem-A(r)}.

\begin{prop}[Volumes]\label{prop-vol-sequence}
For $\Sph^2\times_{f_j}\Sph^1$ as defined in Example \ref{ex-sequence}, 
the volumes are uniformly bounded from above:
\be
\vol_j(\Sph^2\times_{g_j} \Sph^1) \le 4\pi^3\beta.
\ee
\end{prop}

Before we prove Proposition~\ref{prop-vol-sequence}, we prove the following lemma which we will apply both to prove the proposition and to prove other results later in the paper.

\begin{lem}[Areas of Level Sets]\label{lem-A(r)}
Let $A_j(r_0)$ be the area of the level set $r^{-1}(r_0)$ in $\Sph^2\times_{f_j}\Sph^1$ as defined in Example \ref{ex-sequence}.   Then
\be
A_j(r_0)=\int_{\theta=0}^{2\pi} \int_{\varphi=0}^{2\pi} \sin(r_0) f_j(r_0) \,d\varphi\,d\theta=4\pi^2 \sin(r_0)f_j(r_0).
\ee
which increases on $(0,\pi/2)$ and decreases for $r\in (\pi/2,\pi)$ and 
has a maximum value at $r=\pi/2$:
\be
A_j(\pi/2)=4\pi^2 \sin(\pi/2)f_j(\pi/2)=4\pi^2\beta.
\ee
\end{lem}

\begin{proof}
Note that 
\be
f_j'(r)= -\frac{2\cos r\sin r}{\sin^2 r+a_j}.
\ee
Taking the derivative of $A_j(r)$ with respect to $r$, we get
\be
\tfrac{d}{dr}A(r)=4\pi^2(\cos r f_j(r) +\sin r f_j'(r))=4\pi^2\cos r \left( f_j(r)-\frac{2\sin^2 r}{\sin^2 r+a_j}\right)
\ee
Since $f_j(r)> f_j(\pi/2)\geq 2$, the only critical point in $(0,\pi)$ occurs at $r=\pi/2$ where
\be
A_j(\pi/2)=4\pi^2 \sin(\pi/2)f_j(\pi/2) =4\pi^2 1 (\ln(1)+\beta)=4\pi^2\beta.
\ee
Note that
\be \label{SecondDerivativeArea}
\tfrac{d^2}{dr^2}A(r)=4\pi^2(-\sin r f_j(r)+2 \cos rf_j'(r) +\sin r f_j''(r)).
\ee
Recall that by Proposition~\ref{prop-scalar-seq}, the scalar curvature nonnegative condition is equivalent to
\be
\sin r f_j''(r)+\cos r f_j'(r)=\sin r\left(f_j''(r) +\frac{\cos r}{\sin r} f_j'(r)\right)\leq \sin r f_j(r).
\ee 
In particular, using (\ref{SecondDerivativeArea}), we have 
\be
\tfrac{d^2}{dr^2} A_j(r) \leq \cos r f_j'(r)=-\frac{2\cos^2 r\sin r}{\sin^2 r+a_j}\leq 0.
\ee
By the reflective symmetry of $f_j(r)=f_j(\pi-r)$, and the fact that 
\be
\lim_{r\to 0} A_j(r)= \lim_{r\to 0} 4\pi^2 \sin(r)f_j(r) = 0 (\ln((1+a_j)/a_j)+\beta)=0
\ee
we see that $A(r)$ increases on $(0,\pi/2)$,
has a maximum value at $r=\pi/2$ which 
is $\ge 4\pi\beta>4\pi$ and then decreases for $r\in (\pi/2,\pi)$.   
\end{proof}

The proof of Proposition~\ref{prop-vol-sequence} is now easy to see as follows:

\begin{proof}
By the coarea formula,
\be
\vol_j(\Sph^2\times_{g_j} \Sph^1)=\int_{r=0}^\pi\ A_j(r) \, dr
\ee
where $A_j(r)$ is as defined in Lemma~\ref{lem-A(r)}.   Thus by this lemma
\be
\vol_j(\Sph^2\times_{g_j} \Sph^1) \leq  \pi A_j(\pi/2) = 4\pi^3\beta.
\ee
\end{proof}

\subsection{Proving Fibres Stretch to Infinite Length}
\label{subsect-fibres}

Recall that in Lemma~\ref{lem-warp-unbounded}, we proved our warping functions, $f_j$, are unbounded:
\be
\lim_{j\to \infty} f_j(0) \to \infty \textrm{ and } \lim_{j\to \infty} f_j(\pi) \to \infty.
\ee
Here we prove in Lemma~\ref{lem-fibres} that the ${\mathbb S}^1$ fibres over the poles at $r=0$ and $r=\pi$ have lengths increasing to infinity. 

We begin with a review of the standard definition of length of a curve in a Riemannian manifold applied to our sequence:

\begin{defn}[Length of a Curve]\label{DefinitionLengthCurve}
Given a piecewise differentiable curve $c: [a,b]\to \Sph^2\times \Sph^1$, we define the length of $c$ as measured by the metric $g_j$ of Example \ref{ex-sequence}as
\[L_j(c)=\int_{a}^b|\dot{c}(t)|_{g_j}\,dt,\]
where $\dot{c}$ denotes the derivative with respect the $t$.
\end{defn}

\begin{lem}\label{lem-fibres}
The $\Sph^1$ fibres above the poles in Example~\ref{ex-sequence} can be
parametrized as
\be
c_0(t)=(r(t), \theta(t), \phi(t))=(0, 0, t) \textrm{ and }
c_\pi(t)=(r(t), \theta(t), \phi(t))=(\pi, 0, t)
\ee
where $t\in [0, 2\pi]$ and
have
lengths, $L_j(c_0)=2\pi f_j(0)$ and $L_j(c_\pi)=2\pi f_j(\pi)$,
which diverge to infinity as $j \to \infty$.
\end{lem}

\begin{proof}
First note that 
\be
|\dot{c}_0(t)|^2_{g_j}=f^2_j(0) \textrm{ and } |\dot{c}_\pi(t)|^2_{g_j}=f^2_j(\pi).
\ee
By Definition~\ref{DefinitionLengthCurve}
\be
L_j(c_0)=\int_{0}^{2\pi}f_j(0) \, dt=2\pi f_j(0)
\textrm{ and } L_j(c_\pi)=\int_{0}^{2\pi}f_j(\pi) \, dt=2\pi f_j(\pi)
\ee
as claimed.
By Lemma~\ref{lem-warp-unbounded}, these lengths diverge to infinity.
\end{proof}

\subsection{Proving Uniform Distance and Diameter Bounds for the Sequence:}
\label{subsect-unif-dist}

Here we prove the diameter of the sequence of Riemannian manifolds in
Example \ref{ex-sequence} is uniformly bounded from above [Proposition~\ref{UniformBoundDiameter}].  We also prove some bounds on the
distances between points in Lemma~\ref{LemDistanceFunctionsInc} and
Lemma~\ref{LemmaBoundDistance} that we will apply later to prove convergence.

We define the distance function on $(\Sph^2\times_j\Sph^1, g_j)$
of Example \ref{ex-sequence} in the usual way as follows:

\begin{defn}[Distance Function]\label{DefinitionDistanceFunctionSequence}
	For any $p,\ q\in \Sph^2\times\Sph^1$, define
	\be
	d_{j}(p,q)=\inf\{L_j(c)\},
	\ee
	where the infimum is taken over piecewise differentiable curve connecting $p$ and $q$ and $L_j$ is defined as in Definition~\ref{DefinitionLengthCurve}.
\end{defn}

\begin{lem}[Distance Functions are Increasing]\label{LemDistanceFunctionsInc}
 The sequence of distance function $d_j$ as defined in Definition \ref{DefinitionDistanceFunctionSequence} 
 satisfies the inequality:
 	\be
 	d_j(p,q)\leq d_{j+1}(p,q) \qquad \forall j \in {\mathbb N}.
 	\ee
\end{lem}

\begin{proof} 
By Lemma~\ref{lem-warp-inc} we know $g_j \le g_{j+1}$.
	 As a consequence, for any piecewise $C^1$ curve $c$, we have
	$L_j(c)\leq L_{j+1}(c).$
	
	For any $p,\ q\in \Sph^2\times\Sph^1$, since $(\Sph^2\times\Sph^1,g_{j+1})$ is a compact Riemannian manifold with a smooth metric, there exists a minimizing geodesic $c$ connecting $p$ and $q$ such that
	\[d_{j+1}(p,q)=L_{g_{j+1}}(c),\]
	and as a result we have
	\[d_j(p,q)\leq L_j(c)\leq L_{j+1}(c)=d_{j+1}(p,q).\]
\end{proof}

To uniformly bound the diameter, we need to find curves between points whose lengths remain uniformly bounded.  We do this by traveling in the radial direction away from regions where the warping functions are too large.   In the next lemma we show that we can estimate the distance between points involving only the
warping function at one end point.

\begin{lem}[Uniform Bound for Distance Functions]\label{LemmaBoundDistance}
Given any $p_1=(r_1,\theta_1,\varphi_1)$ and $p_2=(r_2,\theta_2,\varphi_2)$ in $\Sph^2\times_{f_j}\Sph^1$
as defined in Example \ref{ex-sequence}, 
	\be
	d_j((r_1,\theta_1,\varphi_1),(r_2,\theta_2,\varphi_2)) \leq 
		|r_1-r_2|+\sin (r_2)\,d_{{\mathbb S}^1}(\theta_1,\theta_2)+f_j(r_2) \,d_{{\mathbb S}^1}(\varphi_1,\varphi_2).
\ee
\end{lem}

\begin{proof}
For any pair of points, $p_1=(r_1,\theta_1,\varphi_1)$ and $p_2=(r_2,\theta_2,\varphi_2)$ in $\Sph^2\times_{f_j}\Sph^1$, 
using triangle inequality we have
\be
\begin{split}
	d_j((r_1,\theta_1,\varphi_1),(r_2,\theta_2,\varphi_2)) \leq & \,\,d_j((r_1,\theta_1,\varphi_1),(r_2,\theta_1,\varphi_1)) \\
	&+d_j((r_2,\theta_1,\varphi_1),(r_2,\theta_2,\varphi_1)) \\
		&+d_j((r_2,\theta_2,\varphi_1),(r_2,\theta_2,\varphi_2)) 
\end{split}
\ee
If we take $c_1(t)=(t,\theta_1,\varphi_1)$ for $t\in [r_1,r_2]$, we have
$|\dot{c}_1(t)|_{g_j}=1$.  So
\be
d_j((r_1,\theta_1,\varphi_1),(r_2,\theta_1,\varphi_1))\le L_j(c_1)= |r_2-r_1|.
\ee
If we take $c_2(t)=(r_2,\theta(t),\varphi_1)$ where $\theta(t)$ is a minimal arc in ${\mathbb S}^1$ from $\theta_1$ to $\theta_2$,
 we have $|\dot{c}_1(t)|_{g_j}=\sin(r_2)$.
So
\be
d_j((r_2,\theta_1,\varphi_1),(r_2,\theta_2,\varphi_1)) \le L_j(c_2)=\sin (r_2)\,d_{{\mathbb S}^1}(\theta_1,\theta_2).
\ee
If we take $c_3(t)=(r_2,\theta_2,\varphi(t))$ where $\varphi(t)$ is a minimal arc in ${\mathbb S}^1$ from $\varphi_1$ to $\varphi_2$,
 we have $|\dot{c}_1(t)|_{g_j}=f_j(r_2)$.  So
\be
d_j((r_2,\theta_2,\varphi_1),(r_2,\theta_2,\varphi_2)) \le L_j(c_3)=f_j(r_2)\,d_{{\mathbb S}^1}(\varphi_1,\varphi_2).
\ee
Combining these four equations completes the proof of the lemma.
\end{proof}

We now apply the above lemma to prove a uniform bound on diameter:

\begin{prop}[Uniform Bound for Diameter]\label{UniformBoundDiameter}
	For $\Sph^2\times_{f_j}\Sph^1$ as defined in Example \ref{ex-sequence}, 
	 we have 
	 \be
	 \diam(\Sph^2\times_{f_j}\Sph^1)\leq (3+2\beta)\pi.
	 \ee
\end{prop}

\begin{proof}
Given any pair of points $p_1, p_2\in \Sph^2\times_{f_j}\Sph^1$ such that $p_2$ has $r(p_2)=\frac{\pi}{2}$,
by Lemma~\ref{LemmaBoundDistance}, we have 
\be \label{diam-1}
d_j(p_1,p_2)\le {\pi}/{2} +\pi +\beta \pi
\ee
because $|r_1-r_2|\leq \frac{\pi}{2}$, $\sin (r_2)\,d_{{\mathbb S}^1}(\theta_1,\theta_2)\leq \pi$ and 
\be
f_j(r_2)\,d_{{\mathbb S}^1}(\varphi_1,\varphi_2)\le (\ln((1+a_j)/(\sin^2(\pi/2)+a_j))+\beta)\pi = \beta \pi. 
\ee
Given any pair of points $q_1, q_2\in \Sph^2\times_{f_j}\Sph^1$, by the triangle inequality,
\be
d_j(q_1,q_2)\le d_j(q_1,p_2)+ d_j(p_2,q_2)
\ee
where $p_2$ has $r(p_2)=\frac{\pi}{2}$.   By (\ref{diam-1}), applied twice, we have
\be
d_j(q_1,q_2)\le ({\pi}/{2} +\pi +\beta \pi)+({\pi}/{2} +\pi +\beta \pi)=(3+2\beta)\pi
\ee
which gives the lemma.
\end{proof}


\newpage
\section{Converging to an Extreme Warped Product Space}
\label{Sect-conv}

We have already shown that the sequence of warping functions, $f_j$, from Example \ref{ex-sequence} are increasing [Lemma~\ref{lem-warp-inc}] and unbounded [Lemma~\ref{lem-warp-unbounded}]   In particular, the fibres above the poles are stretched to infinite length [Lemma~\ref{lem-fibres}].  

In this section we will prove that
 $f_j: {\mathbb S}^2 \to (0,\infty)$ converge smoothly away from the poles to an unbounded  function, $f_\infty: {\mathbb S}^2 \to (0,\infty]$ [Proposition~\ref{prop-warp-conv}].

\begin{figure}[h] 
   \centering
   \includegraphics[width=3in]{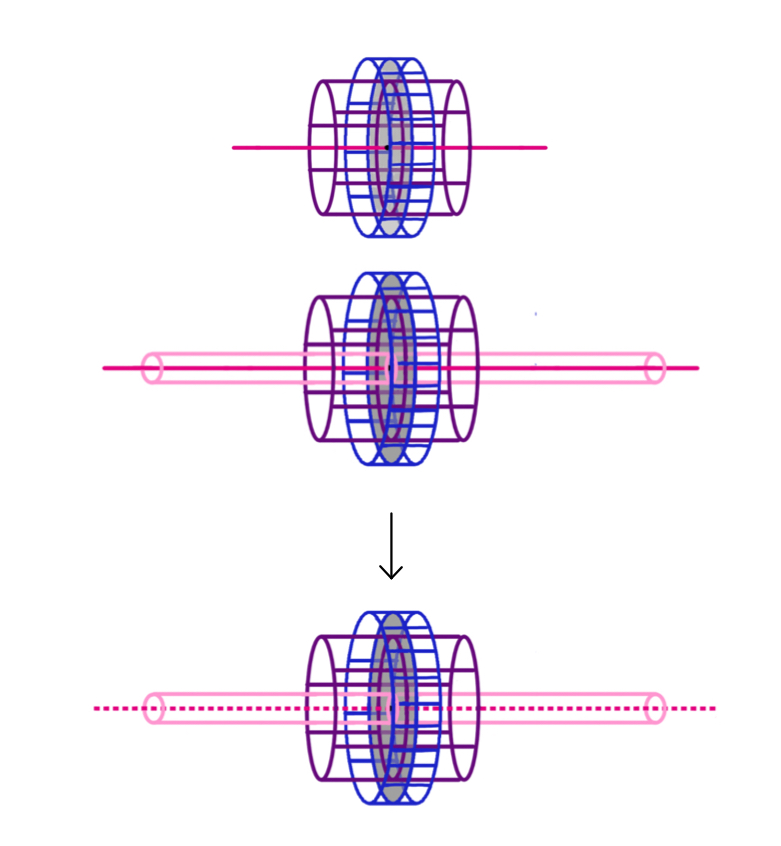} 
   \caption{Here we view fibres above hemispheres in the sequence $\Sph^2\times_{f_j}\Sph^1$ as $j\to \infty$
   so that we can see the warping functions diverging: $f_j(0)\to \infty$. 
   We view our extreme limit space $\Sph^2\times_{f_\infty}\Sph^1$
   with an infinitely stretched ``fibre" above the pole  depicted by a dotted line.  
   }
   \label{fig-STW-2}
\end{figure}

This allows us to define a singular Riemannian manifold, $\Sph^2\times_{f_\infty}\Sph^1$ in Example~\ref{ex-limit} which we call an {\em extreme warped product space} because 
\be
f_\infty(0)=\infty \textrm{ and }f_\infty(\pi)=\infty.
\ee  
Intuitively, this extreme limit space has infinitely stretched "fibres" above the poles.  See Figure~\ref{fig-STW-2}.

\begin{example}\label{ex-limit}
	Consider the extreme warped product space, $\Sph^2\times_{f_\infty}\Sph^1$, 
	with metric tensor
		\be
	g_\infty=g_{\Sph^2}+f_\infty^2(r) g_{\Sph^1}=dr^2+ \sin^2(r) d\theta^2+f_\infty^2(r) d\varphi^2.
	\ee
	written using $(r,\theta)$ coordinates on $\Sph^2$ and $\varphi$ on the $\Sph^1$-fibres
	with the extreme warping function defined by
	\be \label{f-infty}
        f_\infty(r, \theta)=\ln \left(\frac{1}{\sin^2 r}\right)+\beta=-2\ln \sin r+\beta,
        \ee
where $\beta\geq 2$ as in Example \ref{ex-sequence}.   Note that $g_\infty$ is smooth away from the singular set
	 \be
         \label{DefinitionSingularSet}
         S=\{(r,\theta,\varphi)\in\Sph^2\times \Sph^1,\text{ such that } r=0\text{ or }r=\pi\}
         \ee
         so $(\Sph^2\times \Sph^1\backslash S, g_\infty)$ is a smooth open Riemannian manifold.
         See Figure~\ref{fig-STW-2}.
 \end{example}
 
In Subsection~\ref{subsect-smooth-away}, we prove the sequence of metric tensors, $g_j$, in our example converge smoothly on compact sets away from the fibres above the poles to $g_\infty$, in Proposition~\ref{prop-warp-conv}.  Smooth convergence away from singular sets has been studied often before.  See for example the paper of Lakzian-Sormani \cite{Lakzian-Sormani}.

In Subsection~\ref{subsect-Volume}, 
we prove that the volume of the extreme limit is finite in  Lemma~\ref{ExampleVolume}.   So in some sense this is like a cusp singularity, but keep in mind that the singularity is not isolated.  The singular set, $S$, consists of two circles in $\Sph^2\times \Sph^1$.

In Subsection~\ref{subsect-lim-reg} we prove the metric tensor, $g_\infty$, of the extreme limit space 
in Example~\ref{ex-limit} is $W^{1,p}(\Sph^2\times \Sph^1)$ for $p\in [1,2)$ 
but not $H^1_{loc}$ in Proposition ~\ref{metric-W1p}.  In Subsection~\ref{subsect-conv-W1p}, we prove Proposition~\ref{PropConvMetric}
that the metric tensors, $g_j$, of our sequence converge to $g_\infty$ of our extreme limit space in $W^{1,p}(\Sph^2\times \Sph^1)$ as $j\to \infty$ for  $p\in [1,2)$.

\subsection{Smooth Convergence Away from the Singular Set}
\label{subsect-smooth-away}

Here we prove that our sequence of metric tensors, $g_j$, converges to the
metric tensor, $g_\infty$, of the extreme warped space smoothly away from the singular set, $S$, above the poles.   

\begin{prop} \label{prop-warp-conv}
The sequence $g_j$ as defined in Example \ref{ex-sequence} 
converges smoothly away from $S$ to $g_\infty$ of Example~\ref{ex-limit}.   That is, on compact
subsets $K \subset (\Sph^2\times\Sph^1)\setminus S$ we have
\be
|f_j-f_\infty|_{C^\infty(K)}\to 0 \textrm{ and } |g_j-g_\infty|_{g_\infty,C^\infty(K)}\to 0
\ee
\end{prop}

\begin{proof}
For any compact set $K\subset \subset (\Sph^2\times\Sph^1)\setminus S$
there exists $\delta_K$ such that
\be
0<\delta_K\le r(p) \le \pi-\delta_K <\pi \quad \forall p \in K.
\ee
Thus
\be
\sin^2(r(p))+a_j \ge \sin^2(\delta_K)>0 \quad \forall p \in K.
\ee
Since $h_j \to h_\infty$ in $C^\infty(K)$ and $h_j$ uniformly bounded away from $0$ implies
$h_j ^{-1}\to h_\infty^{-1}$ in $C^\infty(K)$ we have 
\be
(\sin^2(r)+a_j)^{-1} \to (\sin^2(r))^{-1} \textrm{ as } a_j \to 0 \textrm{ in } C^\infty(K).
\ee
This implies
\be
\frac{1+a_j}{\sin^2(r)+a_j} +\beta  \to \frac{1}{\sin^2 r}+\beta \textrm{ as } a_j \to 0 \textrm{ in } C^\infty(K).
\ee
which is also bounded uniformly away from $0$:
\be
\frac{1+a_j}{\sin^2(r)+a_j} +\beta  \ge \beta \textrm{ on } K.
\ee
Since
$h_j \to h_\infty$ in $C^\infty(K)$ and $h_j$ uniformly bounded away from $0$ implies
$\ln(h_j )\to \ln(h_\infty)$ in $C^\infty(K)$, we have
\be
\ln\left(\frac{1+a_j}{\sin^2(r)+a_j}\right)+\beta \to \ln \left(\frac{1}{\sin^2 r}\right)+\beta
 \textrm{ as } a_j \to 0 \textrm{ in } C^\infty(K)
\ee
which completes the proof.
\end{proof}

\subsection{The Volume of the Extreme Limit Space:}
\label{subsect-Volume}
 
Here we see that the extreme limit's warping function is Lebesgue and that the volume
of the extreme limit is finite.

\begin{lem}\label{ExampleVolume}   
For any $\beta>0$, we may define $f_\infty: \Sph^2\to \R$ by
\be
f_\infty(r, \theta)=\ln \frac{1}{\sin^2 r}+\beta=-2\ln \sin r+\beta.
\ee
as in (\ref{f-infty}).
Then $f_\infty\in L^1(\Sph^2)$ and
\be
\|f_\infty\|_{L^1(\Sph^2)}=2\pi(2\beta +4-2\ln 4).
\ee
In particular,
\be
\vol_\infty(\Sph^2\times \Sph^1\backslash S)=4\pi^2(2\beta +4-2\ln 4).
\ee
\end{lem}

\begin{proof}
	In the polar coordinate in $\Sph^2$, the metric is
	\be
	g_{\Sph^2}=dr^2+\sin ^2 r d\theta^2.
	\ee
	The volume form is
	\be
	d\vol_{\Sph^2}= \sin(r) dr\wedge d\theta,
	\ee
	and
\be\label{eq-L1}
|f_\infty|_{L^1(\Sph^2)}=\int_{0}^{2\pi}\int_{0}^{\pi} f_\infty(r) \sin r dr d\theta.
\ee
The function $-\ln(\sin r)\sin r$ has anti-derivative
\be
\int -\ln(\sin r)\sin r dr= -\cos r+\ln\left(\cot \frac{r}{2}\right)+\cos r\ln (\sin r)+C.
\ee
After plugging in the bounds, we get
\be
\int_{0}^{\pi} -2\ln(\sin r)\sin r dr=4-2\ln 4.
\ee
As a result we get
\be\label{eq-L1-2}
|f_\infty|_{L^1(\Sph^2)}=2\pi(2\beta +4-2\ln 4).
\ee

The volume of the warped product space can be computed
using its volume form as follows:
\be
\vol_\infty(\Sph^2\times \Sph^1 \backslash S) = \int_{\varphi=0}^{2\pi} \int_{\theta=0}^{2\pi}\int_{r=0}^\pi  \sin(r)d\theta \wedge  f_\infty(r)d\varphi \wedge  dr. 
\ee
Combining this with (\ref{eq-L1}) and (\ref{eq-L1-2}) 
we have
\be
\vol_\infty(\Sph^2\times \Sph^1 \backslash S)= 2\pi |f_\infty|_{L^1(\Sph^2)}
=4\pi^2(2\beta +4-2\ln 4).
\ee
\end{proof}

\subsection{Regularity of the limit metric tensor}
\label{subsect-lim-reg}

Here we prove that the metric tensor, $g_\infty$, of the extreme limit space 
in Example~\ref{ex-limit} is $W^{1,p}$ for all $p\in[1,2)$ but not in $H^1_{loc}$ in Proposition ~\ref{metric-W1p}.   We begin by studying the regularity of the warping function, $f_\infty$, in Lemma~\ref{ExampleRegularity}. 

\begin{lem}[Regularity]\label{ExampleRegularity}
The warping function $f_\infty(r,\theta)=f_\infty(r)$ as defined in Example~\ref{ex-limit} is in $L^{p}(\Sph^2)$ for all $p\in [1,\infty)$, but not in $L^\infty(\Sph^2)$. Also $f_\infty$ is in
$W^{1,p}(\Sph^2)$ for all $p\in[1,2)$, but not in $H^1(\Sph^2)$.
\end{lem}

\begin{proof}
For any $a\geq 0$, $b\geq 0$	 and $p\in[1,\infty)$ we have the elementary inequality
\be
(a+b)^p\leq 2^p(a^p+b^p).
\ee
As a result for any $p\in[1,\infty)$ we have
\be
f_\infty^p\leq 4^p(|\ln(\sin r)|^p+(\beta/2)^p).
\ee
Hence we only need to consider the integrability condition of the function $|\ln(\sin r)|^p$. 

In the polar coordinate in $\Sph^2$ we have
\be
\int_{\Sph^2} |\ln(\sin r)|^p=\int_{0}^{2\pi} \int_{0}^{\pi}|\ln (\sin r)|^p \sin rdrd\theta=4\pi \int_{0}^{\frac{\pi}{2}}|\ln (\sin r)|^p \sin rdr.
\ee
Since the function $|\ln (\sin r)|^p \sin r$ is continuous on the interval $(0,\frac{\pi}{2}]$, we only need to consider the asymptotic behavior when $r\to 0$. Since
\be
\lim_{r\to 0} |\ln (\sin r)|^p \sin r=\lim_{x\to 0} x|\ln x|^p=0,
\ee
for all $p\in [1,\infty)$, we know that $f_\infty\in L^p(\Sph^2)$ for all $p\in [1,\infty)$.

Now we consider the gradient $\nabla f_\infty$. In the polar coordinate in $\Sph^2$, we 
 have
\be
\nabla f_\infty= -2\frac{\cos r}{\sin r}\partial_r,
\ee
and
\be
|\nabla f_\infty|^p=2^p\left|\frac{\cos r}{\sin r}\right|^p.
\ee
When we integrate over $\Sph^2$ we get
\be
\int_{\Sph^2}|\nabla f_\infty|^p= 2^p\int_{0}^{2\pi} \int_{0}^{\pi}\left|\frac{\cos r}{\sin r}\right|^p\sin r dr d\theta=4\pi 2^p\int_0^{\frac{\pi}{2}}\frac{\cos^p r}{\sin ^{p-1}r}dr.
\ee
For any $p\in[0,\infty)$, as $r\to 0$ we have
\be
\lim_{r\to 0}\frac{r^{p-1}}{\sin^{p-1}r}=1,
\ee
where $\frac{1}{r^{p-1}}$ is integrable at $r=0$ only when $p<2$. 
\end{proof}

As a consequence, we will prove the regularity of the metric tensor. But before that we need the following definition:

\begin{defn} \label{defn-W1p}
We define $L^{p}(\Sph^2\times \Sph^1, g_0)$ as the set of all tensors defined almost everywhere on $\Sph^2\times\Sph^1$ such that its $L^p$ norm measured in terms of $g_0$ is finite where $g_0$ is the
isometric product metric 
\be
g_0 = g_{\Sph^2}+ g_{\Sph^1} \textrm{ on } \Sph^2 \times \Sph^1.
\ee
We define $W^{1,p}(\Sph^2\times \Sph^1, g_0)$ as the set of all tensors, $h$,  defined almost everywhere on $\Sph^2\times\Sph^1$ such that both the $L^p$ norm of $h$ and the $L^p$ norm of $\overline{\nabla} h$ measured in terms of $g_0$ are finite
where $\overline{\nabla}$ is the connection corresponding to the metric $g_0$. 
\end{defn}

Now we prove the regularity of the metric tensor $g_\infty$:

\begin{prop} [Regularity of the metric tensor]
\label{metric-W1p}
Choose the product metric $g_0 = g_{\Sph^2}+ g_{\Sph^1}$ on $\Sph^2 \times \Sph^1$ as a smooth background metic, then we have
\be
g_\infty\in W^{1,p}(\Sph^2\times \Sph^1, g_0)
\ee
for all $p\in [1,2)$. 
However, $g_\infty\notin H^1_{loc}(\Sph^2\times \Sph^1, g_0)$.
\end{prop}

\begin{proof}
Using the background metric, we have
\be
\|g_\infty\|^p_{L^p(\Sph^2\times \Sph^1, g_0)}= \int_{\varphi=0}^{2\pi} \int_{\theta=0}^{2\pi}\int_{r=0}^\pi\sin(r) (2+f^4_\infty)^{p/2} dr \, d\theta \, d\varphi.
\ee
From Lemma \ref{ExampleRegularity}, we know that $g_\infty\in L^{p}(\Sph^2\times \Sph^1, g_0)$ for all $p\in [1,\infty)$.

Here we use $\overline{\nabla}$ to denote the connection of the background metric $g_0$. Through direct calculation we have
\be
\overline{\nabla} g_\infty =\overline{\nabla} g_{\Sph^2} +\overline{\nabla} f^2_\infty \otimes g_{\Sph^1} +f^2_\infty \overline{\nabla} g_{\Sph^1}.
\ee
Since $g_0$ is the direct sum metric, we have
\be
\overline{\nabla} g_{\Sph^2}=0,\text{ and }\overline{\nabla} g_{\Sph^1}=0.
\ee
Moreover, since $\overline{\nabla} f^2_\infty= 2f_\infty f'_{\infty}dr$ we have
\be
\overline{\nabla} g_\infty=2f_\infty f'_{\infty}dr\otimes g_{\Sph^1}.
\ee
As  a result, we have
\be
\begin{split}
\|\overline{\nabla} g_\infty\|^p_{L^p(\Sph^2\times \Sph^1, g_0)} &= \int_{\varphi=0}^{2\pi} \int_{\theta=0}^{2\pi}\int_{r=0}^\pi\sin(r) (2f_\infty f'_{\infty})^{p} dr \, d\theta \, d\varphi \\
&=2\pi 2^p \int_{\Sph^2} f_\infty^p (f_\infty^\prime)^p d\vol_{\Sph^2}.
\end{split}
\ee

\textcolor{black}{If $p\in [1,2)$, we can choose $q\in (p,2)$ such that $\frac{q}{p}\in (1,2)$. Choose $t$ such that $\frac{1}{t}+\frac{p}{q}=1$. Note that we have $t\in (2,\infty)$. By H\"older's inequality we have
\be
\begin{split}
\|\overline{\nabla} g_\infty\|^p_{L^p(\Sph^2\times \Sph^1, g_0)} &=2\pi 2^p \int_{\Sph^2} f_\infty^p (f_\infty^\prime)^p d\vol_{\Sph^2}\\
&\leq 2\pi 2^p \left(\int_{\Sph^2} (f_\infty^p)^t d\vol_{\Sph^2}\right)^{1/t}\left( \int_{\Sph^2} (f_\infty^\prime)^{q}d\vol_{\Sph^2} \right)^{p/q}.
\end{split}
\ee
From Lemma \ref{ExampleRegularity} we know that $\|\overline{\nabla} g_\infty\|^p_{L^p(\Sph^2\times \Sph^1, g_0)} < \infty$ and hence $\overline{\nabla} g_\infty\in L^p(\Sph^2\times \Sph^1, g_0)$ for all $p\in[1,2)$.
}

\textcolor{black}{If $p=2$ then we have
\be
\begin{split}
\|\overline{\nabla} g_\infty\|^2_{L^2(\Sph^2\times \Sph^1, g_0)} &= \int_{\varphi=0}^{2\pi} \int_{\theta=0}^{2\pi}\int_{r=0}^\pi\sin(r) (2f_\infty f'_{\infty})^{2} dr \, d\theta \, d\varphi \\
&=16\pi^2  \int_{r=0}^{\pi} \sin r f_\infty^2 (f_\infty^\prime)^2 dr\\
&=16\pi^2  \int_{r=0}^{\pi} \sin r f_\infty^2 \left(2\frac{\cos r}{\sin r}\right)^2 dr\\
&=64\pi^2  \int_{r=0}^{\pi} f_\infty^2\frac{(\cos r)^2}{\sin r} dr.\\
\end{split}
\ee 
Since $ f_\infty^2 \geq \beta^2\geq 4$ and $ \int_{r=0}^{\pi}\frac{(\cos r)^2}{\sin r} dr=\infty $, we conclude that $\|\overline{\nabla} g_\infty\|^2_{L^2(\Sph^2\times \Sph^1, g_0)} =\infty$.
Moreover, since $\Sph^2\times\Sph^1$ is compact, if $g_\infty \in H^{1}_{loc}(\Sph^2\times \Sph^1, g_0)$ then we can conclude that $g_\infty \in W^{1,2}(\Sph^2\times \Sph^1, g_0)$, which is a contradiction. As a result we conclude that $g_\infty\notin H^{1}_{loc}(\Sph^2\times \Sph^1, g_0)$. This finishes the proof.}
\end{proof}

\subsection{Convergence of the Metric Tensor}
\label{subsect-conv-W1p}

Here we prove Proposition~\ref{PropConvMetric}
that the metric tensors, $g_j$, of our sequence converge to $g_\infty$ of our extreme limit space in $W^{1,p}(\Sph^2\times \Sph^1, g_0)$ as $j\to \infty$ for  $p\in [1,2)$.   First we prove convergence of the warping functions in Proposition~\ref{PropConvWarp}.  

\begin{prop}[Convergence of the Warping Function]\label{PropConvWarp}
The sequence $f_j(r,\theta)=f_j(r)$ of Example \ref{ex-sequence} converges to $f_\infty(r,\theta)$ in $L^p(\Sph^2)$ for all $p\in [1,\infty)$ and in
 $W^{1,p}(\Sph^2)$ for all $p\in[1,2)$.
 \end{prop}

\begin{proof}
By Lemma~\ref{lem-warp-inc} and Proposition~\ref{prop-warp-conv}, we know as $j$ increases, $f_j(r,\theta)$ increase to $f_\infty(r,\theta)$
on $r^{-1}(0,\pi)$. 

As functions on $\Sph^2$ we have 
\be
f_j(r,\theta)\le f_\infty(r,\theta) \textrm{ almost everywhere on }\Sph^2.
\ee
\textcolor{black}{Hence we have
\be
|f_j-f_\infty|^p \leq 2^pf_\infty^p, \text{ for all }p\in [1,\infty).
\ee
By Proposition~\ref{prop-warp-conv}, Lemma~\ref{ExampleRegularity} and the Dominated Convergence Theorem, we have $f_j$ converges to $f_\infty$ in $L^p (\Sph^2)$ for all $p\in [1,\infty)$.}

\textcolor{black}{Moreover, since
\be
0\leq |\nabla f_j|=\frac{2|\cos r|\sin r}{\sin^2 r+a_j}\leq \frac{2|\cos r|}{\sin r}=|\nabla f_\infty|.
\ee
Hence for $p\in [1,2)$ we have
\be
|\nabla f_j-\nabla f_\infty|^p \leq 2^p |\nabla f_\infty|^p
\ee
Using Proposition~\ref{prop-warp-conv}, Lemma~\ref{ExampleRegularity} and the Dominated Convergence Theorem, we know that for $p\in[1,2)$
\be \int_{\Sph^2}|\nabla f_j-\nabla f_\infty|^p  \to 0 \text{ as }j\to \infty.
\ee
This finishes the proof.} 
 \end{proof}

\begin{rmrk}
Note that  that $f_j$ does not converge in $W^{1,2}(\Sph^2)$ to $f_\infty$
because we showed $f_\infty \notin W^{1,2}(\Sph^2)$ in Proposition~\ref{ExampleRegularity}.
\end{rmrk}

As a consequence of Proposition \ref{PropConvWarp}, we prove the following 
proposition concerning the metric tensors $g_j$ of Example~\ref{ex-sequence}
and $g_\infty$ of Example~\ref{ex-limit}:

\begin{prop} \label{PropConvMetric}
The metric tensors, $g_j$, converge to $g_\infty$ in $W^{1,p}(\Sph^2\times \Sph^1, g_0)$ as $j\to \infty$ for $p\in [1,2)$.
\end{prop}

Recall the definition of the $W^{1,p}$ norm in Definition~\ref{defn-W1p}.

\begin{proof}
Since
\be
g_j-g_\infty= (f_j^2-f_\infty^2)g_{\Sph^1},
\ee
using the background metric $g_0$ we have
\be
\begin{split}
\|g_j-g_\infty\|^p_{L^p(\Sph^2\times \Sph^1, g_0)}&= 2\pi\int_{\Sph^2} |f_j^2-f_\infty^2|^{p} \\
&= 2\pi\int_{\Sph^2} |f_j-f_\infty|^{p} |f_j+f_\infty|^{p}\\
&\leq  2\pi\int_{\Sph^2} |f_j-f_\infty|^{p} |2f_\infty|^{p}.\\
\end{split}
\ee
Now by Proposition \ref{PropConvWarp} and the H\"older's inequality, we get $g_j$ converges to $g_\infty$ in $L^{p}(\Sph^2\times \Sph^1, g_0)$ as $j\to \infty$ for $p\in [1,\infty)$.

Moreover,
\be
\overline{\nabla}g_j- \overline{\nabla}g_\infty=(2f_j f_j'-2f_\infty f_\infty')dr\otimes g_{\Sph^1}.
\ee
using the background metric $g_0$ we have
\be
\begin{split}
\|\overline{\nabla}g_j-\overline{\nabla}g_\infty\|^p_{L^p(\Sph^2\times \Sph^1, g_0)}&= 2\pi\int_{\Sph^2} 2^p |f_j f_j'-f_\infty f_\infty'|^{p} \\
&= 2^{2p+1}\pi\int_{\Sph^2}\left( |f_j f_j'-f_j f_\infty'|^{p}+|f_j f_\infty'-f_\infty f_\infty'|^{p} \right)\\
&\leq 2^{2p+1}\pi\int_{\Sph^2}\left( |f_j|^p| f_j'- f_\infty'|^{p}+| f_\infty'|^{p} |f_j -f_\infty|^{p} \right)
\end{split}
\ee
where in the second step we used the elementary inequality 
\be
(a+b)^p\leq 2^{p-1} (a^p+b^p)\text{ for all }a,\ b\in [0,\infty).
\ee
Now the conclusion follows from Proposition \ref{PropConvWarp} and the H\"older's inequality. This finishes the proof.
\end{proof}

\section{Curvature of the Extreme Warped Product Space}
\label{Sect-curv}

An essential part of Conjecture~\ref{Scalar-Compactness} is to explore generalized notions of nonnegative scalar curvature on spaces which are the limits of sequences of smooth manifolds with nonnegative scalar curvature.  Thus, in this section, we study the curvature of our limit space.

In Subsection~\ref{subsect-Ricci} and Subsection~\ref{subsect-scalar-away} we study the Ricci and scalar curvature of the smooth part of the extreme warped space of Example~\ref{ex-limit}.   We show the Ricci curvature has no lower bound in Lemma~\ref{LemmaRicci}.   Then applying the smooth convergence away from the singular set, we see that our sequence in Example~\ref{ex-sequence} has no uniform lower bound on Ricci curvature in Proposition~\ref{PropRicci}.
In Remark~\ref{rmrk-not-RCD} we point out that our space thus does not fit in the
classes of spaces with generalized notions of lower bounds on their Ricci curvature.
In Remark~\ref{rmrk-Ricci-flow} we discuss how one might use our calculation to explore
a Ricci flow emerging from our space and check if it has nonnegative scalar curvature in the sense of Burkhardt-Guim \cite{B-G-GAFA}.  

In Subsection~\ref{subsect-scalar-away} we prove that the scalar curvature of
our extreme warped product space is nonnegative where the metric tensor is smooth [Proposition~\ref{LimitScalarCurvature}].

In Subsection~\ref{subsect-distr-scalar} we prove the limit space has nonnegative distributional scalar curvature in the sense of Lee-LeFloch [{defn-Lee-LeFloch}] [Theorem~\ref{thm-distr-scalar}].  Lee-LeFloch defined this notion of
nonnegative scalar curvature in \cite{Lee-LeFloch} building upon the work of LeFloch-Mardare \cite{LM07}.   In Remarks~\ref{rmrk-Lee-LeFloch-original}-\ref{rmrk-Lee-LeFloch}, we discuss how the metric tensors studied by Lee and LeFloch had stronger regularity than
the regularity of $g_\infty$ but nevertheless we are able to prove the distributional scalar curvature is well defined in our setting.

In Subsection~\ref{subsect-total} we  define
total distributional scalar curvature on $\Sph^2\times \Sph^1$ [Definition \ref{defn-distr-total-scalar-curvature}] and compute it for our extreme limit space in Lemma~\ref{lem-distr-total-scalar-curvature}.  We prove that this is the limit of the total scalar curvatures
of our sequence in Proposition~\ref{prop-continuity-of-total-scalar-curvature}.   We also 
compute the total scalar curvature away from the smooth set in Lemma~\ref{lem-total-scalar-curvature-on-smooth-part} and discuss in Remark~\ref{rmrk-total-scalar} what
the total scalar curvature of the singular set must be. 

\subsection{The Ricci Curvature is Unbounded}
\label{subsect-Ricci}

Limits of sequences of manifolds with uniform lower bounds on their Ricci curvature have been studied extensively over the past few decades.   In this section we prove that the sequence in Example~\ref{ex-sequence} has no uniform lower Ricci curvature bound and thus that theory may not be applied to study the sequence.  We complete this proof by first calculating the Ricci curvature of Example~\ref{ex-limit} away from the singular set in Lemma~\ref{LemmaRicci} and then applying the smooth convergence away from the singular set in Proposition~\ref{PropRicci}.
 
\begin{lem}[Ricci Curvature]\label{LemmaRicci}
The Ricci curvature tensor $Ric$ of the warped product $((\Sph^2\times\Sph^1)\backslash S,g_\infty)$ as defined in Example \ref{ex-limit} is not bounded from below.   In particular it satisfies
\begin{enumerate}
	\item $Ric(X,Y)=g_{\Sph^2}(X,Y)-\frac{1}{f} D^2 f_\infty(X,Y)$,
	\item $Ric(X,V)=0$,
	\item $R(V,V)=-\frac{\Delta f_\infty}{f_\infty}g(V,V)=-\frac{\Delta f_\infty}{f_\infty} (f^2g_{\Sph^1}(V,V))$
\end{enumerate}
where $X,\ Y$ to denote vectors tangent to $\Sph^2$, and we use $V$ to denote vectors tangent to $\Sph^1$. $D^2 f_\infty$ denotes the Hessian of $f_\infty$ in the standard $\Sph^2$.  
\end{lem}

\begin{proof}
We can determine the Ricci curvature tensor using Proposition 2 from \cite{KK-compact}.
Note that the Ricci curvature in $\Sph^2$ is the identity metric $g_{\Sph^2}$, and we use $\Delta$ to denote trace of the Hessian, which is different from the notation used in \cite{KK-compact}.
	
	In the polar coordinate we can calculate
\be
D^2f_\infty=\left(\begin{array}{cc}
	2\csc ^2r & 0\\
	0 & -2\cos^2r\\
\end{array}\right).
\ee
As a result, in the orthonormal frame $\{\partial_r,\ \frac{\partial_\theta}{\sin r},\ \frac{\partial_\varphi}{f_\infty}\}$ we can calculate the Ricci tensor to be
\be
Ric=\left(\begin{array}{ccc}
	1-\frac{\csc^2 r}{1-\ln (\sin r)} &0 &0\\
	0 & 1+\frac{\cot^2 r}{1-\ln(\sin r)} &0\\
	0 &0 & -\frac{1}{1-\ln (\sin r)}
\end{array}
\right)
\ee
Note that we have
\be
\begin{split}
\lim_{r\to 0} Ric(\partial_r, \partial_r) &= 1-\lim_{r\to 0} \frac{\csc^2 r}{1-\ln \sin r} \\
&= 1- \lim_{r\to 0}\frac{1}{\sin^2 r(1-\ln \sin r)} \\
&=1-\lim_{x\to 0} \frac{1}{x^2 (1-\ln x)}\\
&=-\infty
\end{split}
\ee
and as a result, the Ricci curvature is not bounded below.
\end{proof}

\begin{prop}[Ricci Curvature of the Sequence]\label{PropRicci}
The Ricci curvature tensor $Ric_j$ of the warped product $((\Sph^2\times\Sph^1),g_j)$ as defined in Example \ref{ex-sequence} is not bounded from below.   That is, there does not exists $H\in {\mathbb R}$ such that for all $j\in {\mathbb N}$ we have
\be\label{eqPropRicci}
\textrm{$Ric_j(v,v)\ge H g_j(v,v)$ for all tangent vectors $v$ at  $p\in (\Sph^2\times\Sph^1)$. }
\ee 
\end{prop}

\begin{proof}
Assume on the contrary that the sequence does satisfy (\ref{eqPropRicci}).   
By Property~\ref{prop-warp-conv}, we have $C^2$ convergence away from the singular set $S$ defined in Example~\ref{ex-limit},
so we have $Ric_j(v,v) \to Ric_\infty(v,v)$.  Thus
\be\label{eqPropRicci}
\textrm{$Ric_\infty(v,v)\ge H g_\infty(v,v)$ for all tangent vectors $v$ at  $p\in (\Sph^2\times\Sph^1)\backslash S$. }
\ee 
This contradicts Lemma~\ref{LemmaRicci}.
\end{proof}

\begin{rmrk}\label{rmrk-not-RCD}
Note that by Lemma~\ref{LemmaRicci} there is no uniform lower bound on the
Ricci curvature and so our space cannot be a CD space as defined
in work of Lott-Villani \cite{Lott-Villani} and Sturm \cite{Sturm}, nor an RCD space as defined in the work of 
Ambrosio-Gigli-Savare \cite{Ambrosio-Gigli-Savare}.  
\end{rmrk}

\begin{rmrk}\label{rmrk-Ricci-flow}
In work of Bamler \cite{Bamler-C}, Burkhardt-Guim \cite{B-G-GAFA}, and  
Huang-Lee \cite{Huang-Lee-scalar} an idea of generalized nonnegative scalar curvature
has been explored in the low regularity setting using Ricci Flow.  Our extreme warped product has lower regularity than they have explored so far.  It would be interesting to explore if there is a smooth solution to Ricci flow $(\Sph^2\times \Sph^1, g_t)$
with nonnegative scalar curvature
such that in the limit as $t\to 0$ we achieve $W^1,p$ convergence of $g_t$ to our
$g_\infty$ of Example~\ref{ex-limit} for $p\in [1,2)$.  Note that Ricci flow starting from singular Riemannian manifolds has been studied by Lee-Topping 
\cite{Lee-Topping-time-zero},  Simon \cite{Simon-deformation} 
Chu-Lee \cite{Chu-Lee-Ricci}, 
Deruelle-Schulze-Simon  \cite{DSS-regularity}, and others.
\end{rmrk}

\subsection{Scalar Curvature of the Limit Space away from the Singular Set}
\label{subsect-scalar-away}

Away from the singular set, our extreme warped product space has a smooth metric tensor, so we can prove the scalar curvature is nonnegative as we did in Subsection~\ref{subsect-scalar} using (\ref{eq-scal-equiv}).   That is
\be
\Scal_\infty=2-2\frac{\Delta f_\infty}{f_\infty}.
\ee
So the nonnegative scalar curvature condition is equivalent to 
\be\label{eq-equiv-lim-scalar}
f_\infty\geq \Delta f_\infty.
\ee

\begin{prop}[Scalar Curvature]\label{LimitScalarCurvature}
The open manifold $((\Sph^2\times\Sph^1)\backslash S,g_\infty)$ as defined in Example \ref{ex-limit} has nonnegative scalar curvature
\be
\Scal_{g_\infty} = 2 - 2\frac{\Delta f_\infty}{ f_\infty } = 2 - \frac{4}{f_\infty},
\ee
\end{prop}

\begin{proof}
In polar coordinate in $\Sph^2$ we have
\be
\Delta f_\infty=\partial_r^2 f_\infty+\frac{\cos r}{\sin r} \partial_r f_\infty.
\ee
We have
\be
 \partial_r f_\infty= - 2\frac{\cos r}{\sin r},
\ee
and
\be
 \partial_r^2 f_\infty= 2 + 2 \frac{\cos^2 r}{\sin^2 r}.
\ee
As a result 
\be
\Delta f_\infty=2.
\ee
On the other hand, since
\be
f_\infty\geq \beta\geq  2,
\ee
we are done by (\ref{eq-equiv-lim-scalar}).
\end{proof}

\subsection{Nonnegative Distributional Scalar Curvature}
\label{subsect-distr-scalar}

Building upon work of Mardare-LeFloch \cite{LM07}, Dan Lee and Philippe LeFloch defined a notion of distributional scalar curvature for smooth manifolds that
have a metric tensor which is only $L^{\infty}_{loc} \cap W^{1, 2}_{loc}$.   See Definition 2.1 of \cite{Lee-LeFloch}
which we review below in Definition~\ref{defn-Lee-LeFloch}.

Here in Theorem~\ref{thm-distr-scalar}, we prove that the metric $g_\infty=g_{\Sph^2}+f_\infty^2 g_{\Sph^1}$ of Example~\ref{ex-limit}
has nonnegative distributional scalar curvature in the sense of Lee-LeFloch .
In Remarks~\ref{rmrk-Lee-LeFloch-original}-\ref{rmrk-Lee-LeFloch}, we discuss how the metric tensors studied by Lee and LeFloch had stronger regularity than
the regularity of $g_\infty$ but their definition of distributional scalar curvature is 
still defined with our lower regularity.

First we recall Definition 2.1 in the work of Lee-LeFloch \cite{Lee-LeFloch}.   In their paper, they assume $M$ is a smooth manifold endowed with a smooth background metric, $g_0$.  They define the scalar curvature distribution, $\Scal_{g}$, as follows for any
metric tensor, $g$, on $M$, such that  $g$ with $L^{\infty}_{loc} \cap W^{1, 2}_{loc}$ regularity and locally bounded inverse $g^{-1} \in L^{\infty}_{loc}$.

\begin{defn}[Lee-LeFloch]\label{defn-Lee-LeFloch} 
The {\em scalar curvature distribution} $\Scal_{g}$ is defined, 
for every test function $u \in C^{\infty}_{0}(M)$, by
\be\label{eqn-Lee-LeFloch}
\langle \Scal_g, u \rangle := \int_{M} \left( - V \cdot \overline{\nabla} \left(u \frac{d\vol_g}{d \vol_{g_0}}\right)  + F u \frac{d\vol_g}{\,d\mu_0}\right) \,d\mu_0,
\ee
where the dot product is taken using the metric $g_0$, $\overline{\nabla}$ is the Levi-Civita connection of $g_0$, $V$ is a vector field given by
\be\label{defn-V}
 V^k:= g^{ij} \Gamma^k_{ij}-g^{ik}\Gamma^j_{ji},
\ee
$F$ is a function as
\be\label{defn-F}
F:= \Scal_{g_0} - \overline{\nabla}_k g^{ij}\Gamma^{k}_{ij} + \overline{\nabla}_{k} g^{ik}\Gamma^{j}_{ji} + g^{ij}\left( \Gamma^{k}_{kl} \Gamma^{l}_{ij} - \Gamma^{k}_{jl}\Gamma^{l}_{ik} \right)
\ee
and
\be\label{eqn-Lee-LeFloch-Christoffel-Sym}
\Gamma^{k}_{ij} := \frac{1}{2} g^{kl} \left( \overline{\nabla}_{i}g_{jl} + \overline{\nabla}_{j}g_{il} - \overline{\nabla}_{l}g_{ij} \right).
\ee
The Riemannian metric $g$ has {\em nonnegative distributional scalar curvature}, if $\langle \Scal_g, u \rangle \geq 0$ for every nonnegative test function $u$ in the integral in (\ref{eqn-Lee-LeFloch}).
\end{defn}

\begin{rmrk}\label{rmrk-Lee-LeFloch-original}
By the regularity assumption for the Riemannian metric $g$ in the work of Lee-LeFloch \cite{Lee-LeFloch}, one has the regularity 
$\Gamma^{k}_{ij} \in L^{2}_{loc}$, $V \in L^{2}_{loc}, F\in L^{1}_{loc}$, 
and the density of volume measure $d\vol_g$ with respect to $\,d\mu_0$ is 
\be
\tfrac{d\vol_{g}}{\,d\mu_0} \in L^{\infty}_{loc} \cap W^{1, 2}_{loc}.
\ee 
Thus
\be\label{LL-term1}
FirstInt_g=\int_{M} \left( - V \cdot \overline{\nabla} \left(u \frac{d\vol_g}{d \vol_{g_0}}\right) \right) \,d\mu_0
\ee
and
\be\label{LL-term2}
SecondInt_g=\int_{M} \left( F u \frac{d\vol_g}{\,d\mu_0}\right) \,d\mu_0.
\ee
are both finite.  Since their sum is the integral in (\ref{eqn-Lee-LeFloch}),
Lee-LeFloch's definition is well defined for the metric tensors that they study.
\end{rmrk}

\begin{rmrk}\label{rmrk-Lee-LeFloch}
Our extreme metric is less regular than the metrics studied by Lee-LeFloch in \cite{Lee-LeFloch}. Recall that in Proposition \ref{metric-W1p} we showed $g_\infty \in W^{1, p}(\Sph^2 \times \Sph^1, g_0)$ for $1\leq p<2$ but $g_\infty \notin W^{1, 2}_{loc}(\Sph^2 \times \Sph^1, g_0)$.
In fact we will show that neither of the integrals in (\ref{LL-term1}) - (\ref{LL-term2}) are finite in Lemma~\ref{lem-extreme-Lee-LeFloch-divergence} below.
However in our Theorem~\ref{thm-distr-scalar} below, we show we have enough regularity to apply Definition~\ref{defn-Lee-LeFloch} because we prove that the integral in (\ref{eqn-Lee-LeFloch}) is finite by summing the integrands first before integrating.
\end{rmrk}

\begin{thm}\label{thm-distr-scalar}
The metric $g_\infty$ has nonnegative distributional scalar curvature on $\Sph^2\times \Sph^1$  in the sense of Lee-LeFloch as in Definition \ref{defn-Lee-LeFloch}.
In particular, (\ref{eqn-Lee-LeFloch}) is finite and nonnegative for any nonnegative
test function,  $u\in C^{\infty}_{0}(\Sph^2 \times \Sph^1)$.
\end{thm}

The proof of Theorem~\ref{thm-distr-scalar} consists of some straightforward but technical calculations. For the convenience of readers, we provide some details of the calculations in the following lemmas, and some of which will be also used in Subsection~\ref{subsect-total} to compute the total scalar curvature of the extreme metric $g_\infty$.

As we have throughout the paper,  we use $g_0 = g_{\Sph^2} + g_{\Sph^1}$ as background metric, and use coordinate $\{r, \theta, \varphi\}$ on $\Sph^2 \times \Sph^1$, where $(r, \theta)$ is a polar coordinate on $\Sph^2$ and $\varphi$ is a coordinate on $\Sph^1$. The corresponding local frame of the tangent bundle is $\{\partial_r, \partial_\theta, \partial_\varphi\}$. In this coordinate system, both $g_0$ and $g_\infty$ are diagonal and given as
\be \label{eqn-diagonal-g_infty}
g_0 =\begin{pmatrix}
1 & 0 &0\\
0& \sin^2 r & 0\\
0&0& 1\\
\end{pmatrix}
\textrm{ and }
g_\infty=\begin{pmatrix}
1 & 0 &0\\
0& \sin^2 r & 0\\
0&0& f_{\infty}^2\\
\end{pmatrix}.
\ee

First of all, by the formula of Christoffel symbols: 
\be
\og^{i}_{jk} = \frac{1}{2} (g_0)^{il}\left( \frac{\partial (g_0)_{il}}{\partial x^k} + \frac{\partial (g_0)_{lk}}{\partial x^j} - \frac{\partial (g_0)_{jk}}{\partial x^l}\right),
\ee
one can easily obtain the following lemma:

\begin{lem}\label{lem-background-metric-Christoffel-sym}
The Christoffel symbols of the Levi-Civita $\on$ of the background metric $g_0 = g_{\Sph^2} + g_{\Sph^1}$, in the coordinate $\{r, \theta, \varphi\}$, all vanish except
\be
\og^{r}_{\theta \theta} = - \sin r \cos r,
\ee
and
\be
\og^{\theta}_{r \theta} = \og^{\theta}_{\theta r} = \frac{\cos r}{\sin r}.
\ee
\end{lem}

Then by Lemma \ref{lem-background-metric-Christoffel-sym}, the formula 
\begin{equation}
\on_i (g_\infty)_{jl}  =  \partial_i \left((g_\infty)_{jl} \right) -  \og^{p}_{ij} (g_\infty)_{pl} - \og^{q}_{il} (g_\infty)_{jq},
\end{equation}
and the diagonal expression of $g_\infty$ in (\ref{eqn-diagonal-g_infty}), one can obtain
the following lemma:

\begin{lem}\label{lem-Lee-LeFloch-Christoffel-sym}
For the extreme metric, $g_\infty$, with the background metric, $g_0$, the Christoffel symbols defined by Lee-LeFloch as in (\ref{eqn-Lee-LeFloch-Christoffel-Sym}), in the coordinate $\{r, \theta, \varphi\}$, all vanish except
\be 
\Gamma^r_{\varphi\varphi}=-f_{\infty}\partial_r f_{\infty},
\ee
and
\be
\Gamma^{\varphi}_{r \varphi}=\Gamma^{\varphi}_{\varphi r}=\frac{\partial_r f_\infty}{f_{\infty}}.
\ee
\end{lem}

Note also that

\begin{lem}\label{lem-vol-ratio}
Note that the volume forms are:
\be
d\mu_0=\,dr\wedge \sin(r)\,d\theta \wedge \, d\phi
\ee
and
\be
d\mu_\infty = dr\wedge \sin(r)\,d\theta \wedge f_\infty(r)\,d\phi
\ee
which are both defined almost everywhere.
In particular,
\be
\frac{\,d\mu_\infty}{\,d\mu_0}= f_\infty(r)
\ee
is in $W^{1,p}(\Sph^2 \times \Sph^1, g_0)$ for $p<2$.
\end{lem}

\begin{proof}
The first two claims hold away from $r=0$ and $r=\pi$ by the definition of volume form.
So $d\mu_\infty = f_\infty d\mu_0$ almost everywhere which gives
us the third claim.  The rest follows from Lemma~\ref{ExampleRegularity}.
\end{proof}

Now we are ready to compute the vector field $V$ and the function $F$ defined by Lee-LeFloch as in (\ref{defn-V}) and $(\ref{defn-F})$.

 \begin{lem}\label{lem-V}
For the extreme metric $g_\infty$ with the background metric $g_0$, the vector field $V$ defined in (\ref{defn-V}), in the local frame $\{\partial_r, \partial_\theta, \partial_\varphi\}$, is given by
\be
V=\left( -2\frac{\partial_r f_\infty}{ f_\infty},0,0\right).
\ee
Furthermore
\be
 - V \cdot \overline\nabla \left(u \frac{\,d\mu_\infty}{\,d\mu_0}\right) 
 =  2 \frac{\partial_r f_\infty}{f_\infty(r)}  \partial_r(u f_\infty (r)).
\ee
\end{lem}

\begin{proof}
By plugging the non-vanishing Christoffel symbols in Lemma \ref{lem-Lee-LeFloch-Christoffel-sym} into
\be
V^k:= g^{ij}_\infty \Gamma^k_{ij}-g^{ik}_\infty\Gamma^j_{ji},
\ee
we get
\begin{eqnarray}
V^r & = &  g^{\varphi \varphi}_\infty \Gamma^r_{\varphi\varphi}- g^{rr}_\infty \Gamma^{\varphi}_{\varphi r} \\
     & = & \frac{1}{(f_\infty)^2} (- f_\infty \partial_r f_\infty) - \frac{\partial_r f_\infty}{f_\infty} = -2 \frac{\partial_r f_\infty}{f_\infty}.
\end{eqnarray}
Also
\begin{equation}
V^{\theta} = g^{ij}_\infty \Gamma^{\theta}_{ij} - g^{\theta \theta}_\infty \Gamma^{j}_{j \theta} = 0.
\end{equation}
\begin{equation}
V^{\varphi} = g^{ij}_\infty \Gamma^{\varphi}_{ij} - g^{\varphi \varphi}_\infty \Gamma^{j}_{j \varphi} = 0.
\end{equation}
By Lemma~\ref{lem-vol-ratio}, we now see that,
\begin{eqnarray}
\overline\nabla \left(u \frac{\,d\mu_\infty}{\,d\mu_0}\right)&=&
\overline\nabla
\left(u f_\infty(r) \right)\\
&=& \partial_r (u f_\infty(r)) + \frac{1}{\sin^2 r} \partial_\theta (u f_\infty(r)) + \partial_\varphi (u f_\infty(r)) 
\end{eqnarray}
Thus
\be
 - V \cdot \overline\nabla \left(u \frac{\,d\mu_\infty}{\,d\mu_0}\right) \\
 =  2 \frac{\partial_r f_\infty}{f_\infty(r)}  \partial_r(u f_\infty (r)).
 \ee
\end{proof}

\begin{lem}\label{lem-F}
For the extreme metric $g_\infty$ with the background metric $g_0$, the function $F$ defined in (\ref{defn-F}) is given by
\begin{equation}
F= 2- 2\left(\frac{\partial_r f_\infty}{f_\infty}\right)^2.
\end{equation}
Furthermore,
\be
\left( F u \frac{\,d\mu_\infty}{d \mu_0} \right)= \left( 2- 2\left(\frac{\partial_r f_\infty}{f_\infty}\right)^2 \right) u f_\infty(r)
\ee
\end{lem}

\begin{proof}
First recall that
\be
\on_i g^{jl}_\infty = \partial_i (g^{jl}_\infty) + \og^{j}_{ip} g^{pl}_\infty + \og^{l}_{iq} g^{jq}_\infty,
\ee
and the scalar curvature of standard product metric $g_0$ is $\Scal_{g_0} = 2$. Then by Lemmas \ref{lem-background-metric-Christoffel-sym} and \ref{lem-Lee-LeFloch-Christoffel-sym}, one has 
\begin{eqnarray}
F & := & \Scal_{g_0} - (\overline{\nabla}_k g^{ij}) \Gamma^k_{ij}+(\overline{\nabla}_k  g^{ik})\Gamma^j_{ji}+g^{ij}(\Gamma^k_{kl}\Gamma^l_{ij}-\Gamma^k_{jl}\Gamma^l_{ik}) \\
& = & 2 - \on_{r}g^{\varphi \varphi} \Gamma^{r}_{\varphi \varphi} - \on_{\varphi}g^{r\varphi} \Gamma^{\varphi}_{r \varphi} - \on_{\varphi} g^{\varphi r} \Gamma^{\varphi}_{\varphi r} + \on_k g^{rk}\Gamma^{\varphi}_{\varphi r} \\
&  & + g^{\varphi \varphi}\Gamma^{\varphi}_{\varphi r} \Gamma^{r}_{\varphi \varphi} - g^{\varphi \varphi} \Gamma^{r}_{\varphi \varphi} \Gamma^{\varphi}_{r \varphi} - g^{rr} \Gamma^{\varphi}_{r \varphi } \Gamma^{\varphi}_{r \varphi} - g^{\varphi \varphi} \Gamma^{\varphi}_{\varphi r} \Gamma^{r}_{\varphi \varphi}\\
& = & 2 - \left( \partial_r (g^{\varphi \varphi}) + \og^{\varphi}_{r \varphi} g^{\varphi \varphi} + \og^{\varphi}_{r \varphi} g^{\varphi \varphi}\right) \Gamma^{r}_{\varphi \varphi}  \\
&  &  - 2 \left( \partial_{\varphi}(g^{r \varphi}) + \og^{r}_{\varphi \varphi} g^{\varphi \varphi} + \og^{\varphi}_{\varphi r} g^{rr} \right) \Gamma^{\varphi}_{r \varphi}\\
& & + \left( \partial_r(g^{rr}) + \og^{r}_{rr}g^{rr} + \og^{r}_{rr} g^{rr} \right) \Gamma^{\varphi}_{\varphi r}  \\ 
& & + \left( \partial_\theta (g^{r \theta}) + \og^{r}_{\theta \theta}g^{\theta \theta} + \og^{\theta}_{\theta r} g^{rr} \right) \Gamma^{\varphi}_{\varphi r} \\ 
&  & + \left( \partial_\varphi (g^{r \varphi}) + \og^{r}_{\varphi \varphi} g^{\varphi \varphi} + \og^{\varphi}_{\varphi r} g^{rr} \right) \Gamma^{\varphi}_{\varphi r} \\
& &  - g^{rr} \Gamma^{\varphi}_{ r \varphi } \Gamma^{\varphi}_{r \varphi} - g^{\varphi \varphi} \Gamma^{\varphi}_{\varphi r} \Gamma^{r}_{\varphi \varphi} \\
& = & 2 - \left( (-2) \frac{\partial_r f_\infty}{(f_\infty)^3} + 0 + 0 \right) \cdot \left( - f_\infty \partial_r f_\infty \right) - 2( 0 + 0 + 0) \Gamma^{\varphi}_{\varphi r} \\
& & + (0 + 0 + 0) \Gamma^{\varphi}_{\varphi r} + \left( 0 + (-\sin r \cos r) \frac{1}{\sin^2 r} + \frac{\cos r}{\sin r} \right) \Gamma^{\varphi}_{\varphi r} \\ 
& & + (0 + 0 + 0) \Gamma^{\varphi}_{\varphi r} \\
& & - \left( \frac{\partial_r f_\infty}{ f_\infty } \right)^2 - \frac{1}{(f_\infty)^2} \frac{\partial_r f_\infty}{f_\infty} (-f_\infty \partial_r f_\infty) \\
& = & 2 - 2 \left( \frac{\partial_r f_\infty}{f_\infty} \right)^2.
\end{eqnarray}
We immediately obtain our second claim by applying Lemma~\ref{lem-vol-ratio}.
\end{proof}
 
\begin{lem}\label{lem-calculation-Lee-LeFloch}
For the extreme $g_\infty = g_{\Sph^2} + f^2_\infty g_{\Sph^1}$, the scalar curvature distribution $\Scal_{g_\infty}$ defined in Definition \ref{defn-Lee-LeFloch} can be expressed, for every test function $u \in C^{\infty}(\Sph^2 \times \Sph^1)$, as the integral
\be
\langle \Scal_{g_\infty}, u \rangle =  \int^{\pi}_{0} -4 \cos r \partial_{r}\bar{u}(r)  +2\bar{u}(r)f_\infty(r) \sin r dr.
\ee
where
\be\label{eq-bar-u}
\bar{u}(r) = \int^{\pi}_{0} \int^{\pi}_{0} u(r, \theta, \varphi) d\theta d\varphi.
\ee
\end{lem}

\begin{proof}   
We choose the product metric $g_0 = g_{\Sph^2}+ g_{\Sph^1}$ on $\Sph^2 \times \Sph^1$ as a smooth background metic. Combining Lemma~\ref{lem-V} and Lemma~\ref{lem-F}
with (\ref{eqn-Lee-LeFloch}) in the definition of distributional scalar curvature, one has
\begin{eqnarray*}
 \langle \Scal_{g_\infty}, u \rangle&=& 
  \int_{\Sph^2\times\Sph^1} \left( - V \cdot \overline\nabla \left(u \frac{\,d\mu_\infty}{\,d\mu_0}\right) + F u \frac{\,d\mu_\infty}{d \vol_{g_0}}\right)\,d\mu_0 \\
& = & \int^{\pi}_{0}  \int^{2\pi}_{0}  \int^{2\pi}_{0} \Bigg[ \,2\,\frac{\partial_r f_\infty(r)}{f_\infty(r)}\,\partial_r \left(u(r, \theta, \varphi)\, f_\infty(r)\right) \\
&  &   \qquad \qquad  \qquad+ \left(2-2\frac{(\partial_{r} f_\infty(r)^2)}{f_\infty^2(r)}\right) \,u(r, \theta, \varphi)\, f_\infty (r)\,\Bigg] \sin r \,d\phi d\theta dr\\
& = & \int^{\pi}_{0} \Bigg[ \,2 \,\frac{\partial_r f_\infty(r)}{f_\infty(r)} \partial_{r} (\bar{u}(r) f_\infty (r)) \\
 &&\qquad\qquad  + \left(2-2\frac{(\partial_{r} f_\infty(r)^2)}{f_\infty^2(r)}\right) \bar{u}(r) f_\infty(r) \Bigg] \sin r dr.
\end{eqnarray*}
Then 
\begin{eqnarray*}
 \langle \Scal_{g_\infty}, u \rangle 
& = &  \int^{\pi}_{0} \left( 2 \frac{\partial_r f_\infty (r)}{f_\infty (r)} \partial_{r} (\bar{u}(r) f_\infty (r)) 
   + \left(2-2\frac{(\partial_{r} f_\infty(r)^2)}{f_\infty^2(r)}\right) \bar{u}(r) f_\infty(r) \right) \sin r dr \\
& = & \int^{\pi}_{0} \Bigg(2 \partial_{r}f_\infty(r) \partial_{r}\bar{u}(r) + 2\frac{(\partial_{r}f_\infty(r))^{2}}{f_\infty(r)}\bar{u}(r) \\
& & \qquad  \qquad\qquad +2\bar{u}(r)f_\infty(r) - 2\frac{(\partial_r f_\infty(r))^2}{f_\infty(r)}\bar{u}(r) \Bigg)\sin r dr \\
& = & \int^{\pi}_{0} 2 \partial_{r}f_\infty(r) \sin r \partial_{r}\bar{u}(r)  +2\bar{u}(r)f_\infty(r) \sin r dr
\end{eqnarray*}
\end{proof}

We now apply these lemmas to prove Theorem~\ref{thm-distr-scalar}:

\begin{proof}
By doing integration by parts for the integral expression of $\langle \Scal_{g_\infty}, u \rangle $ in Lemma \ref{lem-calculation-Lee-LeFloch}, one obtains
\begin{eqnarray*}
 \langle \Scal_{g_\infty}, u \rangle&=&  \int^{\pi}_{0} -4 \cos r \partial_{r}\bar{u}(r)  +2\bar{u}(r)f_\infty(r) \sin r dr \\
&\quad = & -4\cos r \bar{u}(r)|^{\pi}_{0} + \int^{\pi}_{0} -4\sin r \bar{u}(r) + 2\bar{u}(r) f_\infty(r) \sin r dr  \\
& \quad= & 4\bar{u}(0) + 4\bar{u}(\pi)  + \int^{\pi}_{0} (2\beta -4) \sin r \bar{u}(r) + 4 (-\ln \sin r) \sin r \bar{u}(r) dr, \label{RuInnerProduct}
\end{eqnarray*}
since $f_\infty(r)=-2\ln (\sin r) + \beta$. This is nonnegative and finite provided that $u\geq 0$ on
$\Sph^2\times \Sph^1$, because $\beta \geq 2$,
\begin{equation}
4(-\ln \sin r)(\sin r) \bar{u}(r) \geq 0 \textrm{ for $r\in [0, \pi]$,}
\end{equation} 
 and
$$
\lim_{r\to 0} (-\ln(\sin(r))\sin(r)
=
\lim_{s \to 0} (-\ln(s)) s
=
\lim_{h \to \infty} h /e^h
=
\lim_{h \to \infty} 1/e^h=0.
$$
\end{proof}

We close this section by proving the divergence of the integrals in (\ref{LL-term1})
and (\ref{LL-term2}) as discussed in Remark~\ref{rmrk-Lee-LeFloch}:

 \begin{lem}\label{lem-extreme-Lee-LeFloch-divergence}
For $g$ being our extreme metric tensor $g_\infty$ and a smooth nonnegative
test function $u$, the integrals in (150)-(151) are given by
\begin{eqnarray}
\quad FirstInt_{g_\infty}
&=& \int_{\Sph^2 \times \Sph^1}  \left( - V \cdot \overline\nabla \left(u \frac{\,d\mu_\infty}{\,d\mu_0}\right) \right) \,d\mu_0 \\
& = & \int^{\pi}_{0} \left( -4 \cos r \partial_r \bar{u}(r) + 8 H(r) \right) dr, \label{eqn-extreme-metric-FirstInt}
\end{eqnarray}
and
\begin{eqnarray}
\quad SecondInt_{g_\infty}&=& \int_{\Sph^2 \times \Sph^1} \left( F u \frac{\,d\mu_\infty}{d \mu_0} \right) \,d\mu_0\\
& = & \int^{\pi}_{0} \left( 2 \bar{u}(r) f_\infty(r) - 8 H(r) \right) dr,\label{eqn-extrem-metric-SecondInt}
\end{eqnarray}
where
\be
H(r)=\left(\frac{\cos^2 r}{\sin r}\right)\left(\frac{\bar{u}(r)}{(-2 \ln \sin r + \beta)}\right)\ge 0
\ee
and where $\bar{u}$ is defined as in (\ref{eq-bar-u}).

Moreover, if a smooth nonnegative test function $u$ is chosen so that at least one of $\bar{u}(0)$ and $\bar{u}(\pi)$ is non-zero, then neither of integrals in (\ref{eqn-extreme-metric-FirstInt}) and (\ref{eqn-extrem-metric-SecondInt}) are finite.
\end{lem}

\begin{proof}
By integrating the formulas in Lemma~\ref{lem-V} and Lemma~\ref{lem-F}, one can easily obtain the integrals in (\ref{eqn-extreme-metric-FirstInt}) and (\ref{eqn-extrem-metric-SecondInt}). 

Clearly, for any smooth test function $u$, the integrals 
\be
\int^{\pi}_{0} -4\cos r \,\partial_r \bar{u}(r)\,dr \textrm{ and }
\int^{\pi}_{0} 2 \bar{u} \,f_\infty(r) \,dr
\ee
 are both finite. 
 If $\bar{u}(0) \neq 0$ or $\bar{u}(\pi) \neq 0$, then we claim 
 \be
 \int_0^\pi H(r) \,dr
 \ee
 is an improper integral with the convergence issue at $r=0$ or $r=\pi$ respectively. Without loss of generality, we assume $\bar{u}(0) >0$. Then by the continuity of $\bar{u}(r)$, there exists $0 < r_0 < \frac{\pi}{4}$ such that $\bar{u}(r) \geq \frac{\bar{u}(0)}{2}$ for all $r \in [0, r_0]$. So
\begin{eqnarray*}
 \int_0^\pi &H(r)& \,dr \,\,\ge\,\, \int^{r_0}_{0} H(r)\,dr\\
 &=& \int^{r_0}_{0} \left(\frac{\cos^2 r}{\sin r}\right)\left(\frac{\bar{u}(r)}{(-2 \ln \sin r + \beta)}\right)\, dr  \\
& = & \lim_{\epsilon \rightarrow 0} \int^{r_0}_{\epsilon} \frac{\cos^2 r}{\sin r} \frac{\bar{u}(r)}{(-2\ln \sin r + \beta)} dr\\
& \geq & \frac{\bar{u}(0)}{2} \cdot  \cos r_0 \cdot \lim_{\epsilon \rightarrow 0} \int^{r_0}_{0} \frac{\cos r}{\sin r} \frac{1}{(-2 \ln \sin r + \beta)} dr \\
& = &  \frac{\bar{u}(0)}{2} \cdot  \cos r_0 \cdot \lim_{\epsilon \rightarrow 0}  \frac{\left[ \ln (-2 \ln \sin \epsilon + \beta)  - \ln(-2 \ln \sin r_0 + \beta) \right]}{2} \\
& = & +\infty.
\end{eqnarray*}
\end{proof}

\subsection{Distributional Total Scalar Curvature}
\label{subsect-total}

Recall the singular set of the extreme metric $g_\infty$ of Example~\ref{ex-limit} is
\be
S = \{ (r, \theta, \varphi)  \in \Sph^2 \times \Sph^1 \mid r=0, \pi \},
\ee
and $g_\infty$ is smooth on the regular part  $\Sph^2 \times \Sph^1 \setminus S$. So in Lemma~\ref{lem-total-scalar-curvature-on-smooth-part}, we compute the total scalar curvature on the regular part by integrating the usual scalar curvature.   However, intuitively there
seems to be an infinite amount of scalar curvature concentrated on the singular set.

In Definition \ref{defn-distr-total-scalar-curvature} we define a notion of total distributional scalar curvature building upon the definition of Lee-LeFloch from the previous subsection.
We compute the total distributional scalar curvature on $\Sph^2\times \Sph^1$ in Lemma~\ref{lem-distr-total-scalar-curvature}.  In Proposition~\ref{prop-continuity-of-total-scalar-curvature}, we prove the total scalar curvatures of the metrics in the sequence defined in Example \ref{ex-sequence} converge to the distributional total scalar curvature of the extreme limit metric defined in Example \ref{ex-limit}.
In Remark~\ref{rmrk-total-scalar} we discuss
how the difference between the total distributional scalar curvature on $\Sph^2 \times \Sph^1$ and the total smooth scalar curvature over $\Sph^2 \times \Sph^1 \setminus S$
can be explained as the contribution of the singular set to the total scalar curvature. 

\begin{lem}\label{lem-total-scalar-curvature-on-smooth-part}
The total scalar curvature of $g_\infty$ on the regular part $\Sph^2 \times \Sph^1 \setminus S$ is given by
\be\label{eqn-total-scalar-curvature-on-smooth-part}
\int_{\Sph^2 \times \Sph^1 \setminus S} \Scal_{g_\infty} \,d\mu_\infty = (2\pi)^2 \int^{\pi}_{0} 2 \sin r f_\infty(r) dr - 8 (2\pi)^2.
\ee
\end{lem}

\begin{proof}
By Proposition~\ref{LimitScalarCurvature},
on the regular part $\Sph^2 \times \Sph^1 \setminus S$, 
\be
\Scal_{g_\infty} = 2 - 2\frac{\Delta f_\infty}{ f_\infty } = 2 - \frac{4}{f_\infty}.
\ee
Thus  
\begin{eqnarray}
 \int_{\Sph^2 \times \Sph^1 \setminus S} &\Scal_{g_\infty} &\,d\mu_\infty= \\
 & = & \int^{2\pi}_{0} d\theta \int^{2\pi}_{0} d \varphi \int^{\pi}_{0} \left( 2 -\frac{4}{f_\infty (r)} \right) \sin r f_\infty(r) dr \\
 & = & (2\pi)^2 \int^{\pi}_{0} (2 \sin r f_\infty(r) - 4 \sin r )dr \\
 & = & (2\pi)^2 \int^{\pi}_{0} 2 \sin r f_\infty(r) dr - (2\pi)^2 4 \int^{\pi}_{0} \sin r dr \\
 & = & (2\pi)^2 \int^{\pi}_{0} 2 \sin r f_\infty(r) dr - 8 (2\pi)^2. 
\end{eqnarray}
\end{proof}

Now in order to give a description for the total scalar curvature on the whole manifold, including singular set, we would like to use the scalar curvature distribution defined by Lee-LeFloch as in Definition \ref{defn-Lee-LeFloch}. 

\begin{defn}[Distributional total scalar curvature]\label{defn-distr-total-scalar-curvature}
The total distributional scalar curvature of metric $g$ is
 $\langle \Scal_{g}, 1 \rangle$, which is obtained by setting the test function $u\equiv 1$ in the integration in (\ref{eqn-Lee-LeFloch}). 
 \end{defn}

\begin{rmrk}\label{rmrk-total-LL}
On the smooth compact Riemanian manifold, the distributional total scalar curvature is the same as the usual total scalar curvature, i.e. the integral of the scalar curvature. In general, on a non-compact manifold $M$, for a Riemannian metric with the regularity assumptions as in Definition \ref{defn-Lee-LeFloch}, the distributional total scalar curvature may diverge. But for our extreme metric $g_\infty$ on $\Sph^2 \times \Sph^1$, it is finite.
\end{rmrk}

\begin{lem}\label{lem-distr-total-scalar-curvature}
The extreme metric $g_\infty$ defined in Example \ref{ex-limit} has finite and positive distributional total scalar curvature as 
\be\label{eqn-distr-total-scalar-curvature}
\langle \Scal_{g_\infty}, 1 \rangle = (2 \pi)^2 \int^{\pi}_{0} 2 \sin r f_{\infty}(r) dr.
\ee
\end{lem}

\begin{proof}
Note that for $u \equiv 1$ one has $\bar{u} =  \int^{2\pi}_{0}\int^{2\pi}_{0} u(r, \theta, \varphi)d\theta d\varphi \equiv (2\pi)^2$. Then by using Lemma \ref{lem-calculation-Lee-LeFloch}, we obtain
\be
\langle \Scal_{g_\infty}, 1 \rangle = (2 \pi)^2 \int^{\pi}_{0} 2 \sin r f_{\infty}(r) dr.
\ee
Then because $f_\infty (r) = -2 \ln \sin r + \beta \geq 0$, and $\sin r \ln \sin r \rightarrow 0$ as $r \rightarrow 0$ or $\pi$, we see that the integral in (\ref{eqn-distr-total-scalar-curvature}) is convergent and nonnegative.
\end{proof}

\begin{rmrk} \label{rmrk-total-scalar}
{\rm
Clearly distributional total scalar curvature given in (\ref{eqn-distr-total-scalar-curvature}) is strictly greater than the total scalar curvature on the regular part given in (\ref{eqn-total-scalar-curvature-on-smooth-part}). The difference between them could be explained as the contribution of the singular set $S$ to the total scalar curvature, which is positive and equal to $8(2\pi)^2$.
}
\end{rmrk}

\begin{prop}[Continuity of distributional total scalar curvature]\label{prop-continuity-of-total-scalar-curvature}
The total scalar curvatures of the metrics in the sequence defined in Example \ref{ex-sequence} converge to the distributional total scalar curvature of the extreme limit metric defined in Example \ref{ex-limit}.
\end{prop}

\begin{proof}
By using the scalar curvature formula for metrics $g_j$ given in Proposition \ref{prop-scalar-seq},
we have
\begin{eqnarray}
 \int_{\Sph^2 \times \Sph^1} \Scal_{g_j} \,d\mu_j&=&
\int^{2\pi}_{0} d\varphi \int_{\Sph^2} \left( 2 -2 \frac{\Delta f_j}{f_j} \right) f_j \,d\mu_{\Sph^2}\\
& = & (2\pi) \int_{\Sph^2} (2 f_j -2\Delta f_j) \,d\mu_{\Sph^2}\\
& = & (2\pi) \int_{\Sph^2} 2 f_j \,d\mu_{\Sph^2}.
\end{eqnarray}
By the $L^1$ convergence of $f_j$ to $f_\infty$ obtained in Proposition \ref{PropConvWarp}, we see that
\begin{eqnarray}
\lim_{j\to \infty} \int_{\Sph^2 \times \Sph^1} \Scal_{g_j} \,d\mu_j&=&
 (2\pi)   \int_{\Sph^2} 2 f_\infty \,d\mu_{\Sph^2} \\
& = & (2\pi)^2 \int^{\pi}_{0} 2  f_\infty(r) \sin(r) \,dr 
\end{eqnarray} 
which is the total distributional scalar curvature by Lemma~\ref{lem-distr-total-scalar-curvature}.
 \end{proof}

\section{Proving the Sequence has a Uniform Lower Bound on MinA:}
\label{sect-MinA}

Here we study the minimal surfaces in our sequence in Example~\ref{ex-sequence}.  We prove the following theorem:

\begin{thm}[Uniform Bound for MinA]\label{thm-MinA}
There
is a uniform positive lower bound, $A_0>0$ such that any closed $g_j$-minimal surface, $\Sigma$, in 
$\Sph^2\times_{f_j}\Sph^1$ has area $\area_j(\Sigma)\ge A_0$.   Thus 
\be
MinA(\Sph^2\times_{f_j}\Sph^1)\ge A_0>0 \quad \forall j \in {\mathbb N}
\ee
\end{thm}

Before we prove Theorem~\ref{thm-MinA}, we prove various lemmas about minimal surfaces in $\Sph^2\times_{f_j}\Sph^1$ in Subsection~\ref{subsect-min-surf}.  We then
apply some of these lemmas to prove the theorem in Subsection~\ref{subsect-pf-MinA}.
We close with further remarks about the properties of minimal surfaces in $\Sph^2\times_{f_j}\Sph^1$ and conjecture that $A_0=4\pi$ in Subsection~\ref{subsect-rmrks}.  

It is worth noting that minimal surfaces and their relationship with nonnegative scalar curvature was initially studied by Schoen-Yau in \cite{Schoen-Yau-minimal}.  Since then, this relationship has been explored by many mathematicians.  See the survey by Zhu
\cite{Richard-Zhu-small}  and the survey by the first author \cite{Sormani-conjectures}.
See in particular the work of Marques-Neves \cite{MN-Duke}, Bray-Brendle-Neves \cite{BBN10}, Gromov-Zhu \cite{Gromov-Zhu-area}, Liokumovich-Maximo \cite{LioMax-Waist},  and Zhu \cite{Zhu-2inNdim}.    In this paper we do not directly apply these results.  In fact we do not use the nonnegative scalar curvature when studying the minimal surfaces.

\subsection{Minimal Surfaces and Mean Convexity}
\label{subsect-min-surf}

Recall that by Definition~\ref{defn-MinA}, 
\be
\MinA (\Sph^2\times_{f_j}\Sph^1)=\min \area_j(\Sigma)
\ee
where $\Sigma$ is a $g_j$-minimal surface, $\Sigma\subset \Sph^2\times_{f_j}\Sph^1$.  More precisely $\Sigma$ is a 
surface without boundary lying in $\Sph^2\times_{f_j}\Sph^1$ which has  $H_j(\Sigma, p,\nu)=0$ at every $p\in \Sigma$
where $H_j(\Sigma, p, \nu)$ is the mean curvature defined using the metric tensor $g_j$ at $p\in \Sigma$ with normal $\nu$.   

\begin{figure}[h] 
   \centering
   \includegraphics[width=4in]{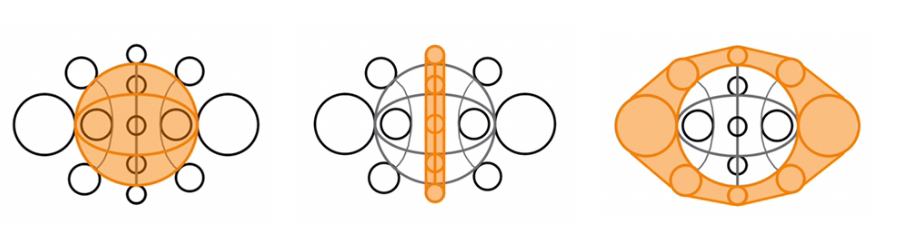} 
   \caption{Three minimal surfaces in $\Sph^2\times_{f_j}\Sph^1$.   On the left is
   $\Sigma_{\varphi=\varphi_0}$ of Lemma~\ref{lem-Sigma-phi}, in the
   middle is $\Sigma_{r=\pi/2}$ of Lemma~\ref{lem-Sigma-r}, and
    on the right is $\Sigma_{\theta=\theta_0}$ of Lemma~\ref{lem-Sigma-theta}.
   }
   \label{fig-STW-3-min}
\end{figure}

In general a surface might be $g_j$-minimal for only one
value of $j\in {\mathbb N}$, however we show there are three classes of  surfaces which
are $g_j$-minimal for all values of $j\in {\mathbb N}$.  These are depicted in Figure~\ref{fig-STW-3-min}
and are proven to be $g_j$-minimal with $\area_j\ge 4\pi$ in Lemmas~\ref{lem-Sigma-phi}-~\ref{lem-Sigma-theta}. 

\begin{lem} \label{lem-Sigma-isom}
If there is a $g_j$-isometry $I: M\to M$, $I_*g_j=g_j$,
that preserves $\Sigma$,
$I(p)=p \,\, \forall p\in \Sigma$, but reverses its orientation,
\be
I^* \nu= -\nu \textrm{ where } \nu \perp \Sigma,
\ee
then $\Sigma$ is a $g_j$-minimal surface.
\end{lem}

\begin{proof} 
The $g_j$ mean curvature of $\Sigma$ at $p\in \Sigma$
depends on a choice of normal, $\nu$, at $p$ as follows:
\be
H_j(\Sigma,p,\nu)=-H_j(\Sigma,p,-\nu).
\ee
Since $I$ is an isometry preserving $\Sigma$ and $I(p)=p$ we have
\be
H_j(\Sigma, p,\nu)=H_j(I(\Sigma), I(p), (I^*\nu))=H_{j}(\Sigma, p,-\nu)
\ee
which implies $H_j=0$.
\end{proof}

\begin{lem} \label{lem-Sigma-phi}
For any $\varphi_0\in {\mathbb S}^1$, the surface, 
\be
\Sigma_{\varphi=\varphi_0}^2=\{\,(r,\theta,\varphi)\,:\, \varphi=\varphi_0\,\} \,\subset\, \Sph^2\times_{f_j}\Sph^1,
\ee
which is a sphere, is a minimal surface of area $4\pi$.  See Figure~\ref{fig-STW-3-min}.
\end{lem}

\begin{proof}
This follows by Lemma~\ref{lem-Sigma-isom} because the reflection across the surface, $I(r,\theta,\varphi)=(r, \theta, (2\varphi_0-\varphi)mod 2\pi)$,
is an isometry:
\be
I_*g_j=I_*(dr^2+\sin^2(r) d\theta^2+ f_j^2(r) d\varphi^2)=dr^2+\sin^2(r) d\theta^2+ f_j^2(r) d\varphi^2.
\ee
The area is computed using its parametrization by $r$ and $\theta$:
\be
\area(\Sigma_{\varphi=\varphi_0}^2)=\int_{r=0}^\pi \int_{\theta=0}^{2\pi} \sin(r)\,d\theta\,dr=\area({\mathbb S}^2)=4\pi.
\ee
\end{proof}

\begin{lem} \label{lem-Sigma-theta}
For any $\theta_0\in {\mathbb S}^1$, the surface, 
\be
\Sigma_{\theta=\theta_0}^2=\{\,(r,\theta,\varphi)\,:\, \theta=\theta_0 \textrm{ or } \theta=(\pi+\theta_0) mod 2\pi\,\} \,\subset\, \Sph^2\times_{f_j}\Sph^1,
\ee
which is a torus above a longitude,
is a minimal surface of area $>4\pi$.  
\end{lem}

\begin{proof}
This follows by Lemma~\ref{lem-Sigma-isom} because the reflection across the surface, $I(r,\theta,\varphi)=(r, (2\theta_0-\theta)mod2\pi, \varphi)$,
is an isometry:
\be
I_*g_j=I_*(dr^2+\sin^2(r) d\theta^2+ f_j^2(r) d\varphi^2)=dr^2+\sin^2(r) d\theta^2+ f_j^2(r) d\varphi^2.
\ee
The area is computed using its parametrization by $r$ and $\varphi$:
\be
\area(\Sigma_{\theta=\theta_0}^2)=\int_{r=0}^\pi \int_{\varphi=0}^{2\pi} f_j(r) \,d\varphi\,dr>2\pi^2 \beta>4\pi.
\ee
\end{proof}

\begin{lem} \label{lem-Sigma-r}
For any $r_0\in (0,\pi)$ the surface, 
\be
\Sigma_{r=r_0}^2=\{\,(r,\theta,\varphi):\, r=r_0\,\}\subset \Sph^2\times_{f_j}\Sph^1
\ee
which is a torus above the latitude $r=r_0$,
is a surface of constant mean curvature with area 
\be
\area_j(\Sigma_{r=r_0})= 4\pi \sin(r_0)f_j(r_0).
\ee
When $r_0=\pi/2$,  it 
is a $g_j$-minimal surface of area $\area_j(\Sigma_{r=\pi/2})>4\pi$ as in Figure~\ref{fig-STW-3-min}. 
When $r_0\neq \pi/2$
the mean curvature is positive with normal chosen pointing away from the poles towards the equator.
\end{lem}

\begin{proof}
The surface, $\Sigma_{r=r_0}$ for $r_0\in (0,\pi)$ is a torus, 
\be
\Sigma_{r=r_0}=\{(r,\theta,\pi)\,:\, \theta\in {\mathbb S}^1, \varphi \in {\mathbb S}^1\}.
\ee
and there are $g_j$-isometries acting transitively on it because $g_j$ is constant in $\theta$ and $\varphi$. Thus it has constant mean curvature.  The area is computed using its parametrization by $\theta$ and $\varphi$:
\be
\area_j(\Sigma_{r=r_0}^2)=A_j(r_0)=\int_{\theta=0}^{2\pi} \int_{\varphi=0}^{2\pi} \sin(r_0) f_j(r_0) \,d\varphi\,dr=4\pi^2 \sin(r_0)f_j(r_0).
\ee
By Lemma~\ref{lem-A(r)} we know $A_j(r)$ increases on $(0,\pi/2)$,
has a maximum value at $r=\pi/2$ which 
is $4\pi\beta>4\pi$ and then decreases for $r\in (\pi/2,\pi)$. 
Thus $\Sigma_{r=r_0}$ is mean convex away from the poles
towards the equator for $r_0\neq \pi/2$ and is an unstable $g_j$-minimal surface when $r_0=\pi/2$.
\end{proof}

The following well known proposition follows from the maximum principal and is easily verified:

\begin{prop} \label{prop-min-K}
If $\Sigma\subset K$ is a minimal hypersurface of codimension one in a Riemannian manifolds and $K$ is a smooth compact region and there is a point $p\in \Sigma \cap \partial K$ 
then $H(\partial K, p, \nu)\le 0$ measured using the outward normal $\nu$.
\end{prop}

The following lemma applies this proposition.

\begin{lem}\label{lem-hit-eq}
If $\Sigma$ is a $g_j$-minimal surface then $\Sigma \cap \Sigma_{r_0=\pi/2} \neq \emptyset$.
\end{lem}

\begin{proof}
Assume on the contrary that $\Sigma \cap \Sigma_{r_0=\pi/2} = \emptyset$.  Without loss of generality,
$\Sigma \subset r^{-1}[0,\pi/2)$.   Let $r_0$ be the maximum value of $r$ of $\Sigma$.   Let $K=r^{-1}[0,r_0]$.
Then by Lemma~\ref{lem-Sigma-r}, 
$H(\partial K, p, \nu)> 0$ measured using the outward normal $\nu$.  This contradicts the maximum principal, Proposition~\ref{prop-min-K}.
\end{proof}

\subsection{Proof of Theorem~\ref{thm-MinA}:}
\label{subsect-pf-MinA}

We will apply the Monotonicity Formula (7.5) from Colding \& Minicozzi's 
textbook \cite{CoMin-text} to prove Theorem~\ref{thm-MinA}, so we state it here:

\begin{thm}[Monotonicity Formula]\label{CM-JJ:thm}
Let $p$ be a point on a smooth minimal surface $\Sigma'$ in a $3$-manifold $(M,g)$. Let $\kappa>0$ be a bound on sectional curvatures $K_M$ on $M$ (as in $|K_M|<\kappa$) and let $i_0>0$ denote a positive lower bound on the injectivity radius on $M$. Then the function  
\be
e^{2\sqrt{\kappa}s} s^{-2} \area(B_{g}(p,s)\,\cap \, \Sigma')
\ee
is non-decreasing
for $s\in \left(0,\min\{i_0, \kappa^{-1/2}, \mathrm{dist}(p, \partial \Sigma')\}\right)$.   Here the restriction that $s<\mathrm{dist}(p, \partial \Sigma')$
is only needed if $\Sigma'$ has a boundary.
\end{thm}

We now prove Theorem~\ref{thm-MinA}:

\begin{proof}
Assume on the contrary, that there exists a subsequence $j\to \infty$ and $g_j$-minimal surfaces $\Sigma_j \subset \Sph^2\times_{f_j}\Sph^1$
such that $\area_j(\Sigma_j)\to 0$.   By Lemma~\ref{lem-hit-eq}, there exists $p_j\in \Sigma_j$ such that $r(p_j)=\pi/2$.   In particular 
\be
\Sigma_j' =\Sigma_j \cap r^{-1}[\pi/4, 3\pi/4] \neq \emptyset
\ee
is a $g_j$-minimal surface which possibly has a boundary such that if $\partial \Sigma'_j\neq \emptyset$ then
\be
\mathrm{dist}_{g_j}(p_j, \partial \Sigma'_j)\ge \min\left\{ |r(p_j)-r(q)|\,:\, q\in \partial \Sigma_j'\right\} =\pi/4.
\ee
Let 
\be
\hat{g}_j=dr^2+\sin^2(r) d\theta^2 + \hat{f}^2_j(r) d\varphi^2
\ee
where $\hat{f}_j(r)=\ln\left(\frac{1+a_j}{S(r)+a_j}\right)+\beta$ where $S(r(x))$ is a smooth positive function on $\Sph^2\times\Sph^1$
such that
\be
S(r)=\sin^2(r) \textrm{ for }r\in (\pi/4, 3\pi/4).
\ee
Thus $\Sigma_j'$ is also a $\hat{g}_j$-minimal surface which possibly has a boundary with
\be
\mathrm{dist}_{\hat{g}_j}(p_j, \partial \Sigma'_j)\ge \pi/4.
\ee
Furthermore, for $s<\pi/4$, $B_{\hat{g}_j}(p_j,s)=B_{g_j}(p_j,s)$ and 
\be\label{area-match}
\area_{\hat{g}_j}\left(\Sigma'_j\cap B_{\hat{g}_j}(p_j,s)\right)=\area_{{g}_j}\left(\Sigma_j\cap B_{{g}_j}(p_j,s)\right) \le \area_j\left(\Sigma_j\right).
\ee
Recall that by Lemma~\ref{prop-warp-conv} that as $j \to \infty$, $f_j(r)$
converges smoothly away from the poles to $f_\infty(r)$.

Thus $\hat{f}_j(r)$ converge smoothly and uniformly to the smooth positive function $\hat{f}_\infty(r)=-\ln(S(r))+\beta$.
So $\hat{g}_j$ converge smoothly and thus have uniform bounds on their sectional curvatures, $\kappa_j\le \kappa_0$, and on their injectivity radii, $i_{0,j}\ge i_0>0$.   In particular. there is an $s_0>0$ such that
\be
s_0< \min\left\{i_{0,j},\kappa_j^{-1/2}, \pi/4\right\} \textrm{ for all } j \in {\mathbb N}.
\ee
We do not have such uniform bounds on the original metric tensors $g_j$.

Applying the Monotonicity Formula recalled in Theorem~\ref{CM-JJ:thm} to $\Sigma'_j$ as $\hat{g}_j$-minimal surfaces, we have for any $s\in (0,s_0)$,
\be
e^{2\sqrt{\kappa_j}s} s^{-2} \area_{\hat{g}_j}\left(B_{\hat{g}_j}(p,s)\,\cap \, \Sigma'_j\right) \ge 4\pi
\ee
because $4\pi$ is the limit as $s\to 0$.   By (\ref{area-match}), we have
\be
\area_j\left(\Sigma_j\right)\ge 4\pi s_0^2 e^{-2\sqrt{\kappa_0}s_0}, 
\ee
which contradicts $\area_j(\Sigma_j) \to 0$.
\end{proof}

\subsection{Remarks and Conjectures}
\label{subsect-rmrks}

Recall that we proved $\Sigma_{\varphi=\varphi_0}$,   $\Sigma_{\theta=\theta_0}$,
and $\Sigma_{r=\pi/2}$ all have area $\ge 4\pi$ in Lemma~\ref{lem-Sigma-phi},
Lemma~\ref{lem-Sigma-theta}, and Lemma~\ref{lem-Sigma-r} respectively.  In addition, we observe that if there is a surface $\Sigma$ whose area is $<4\pi$ then the projection of the surface to the sphere by $(r,\theta,\varphi)\to (r,\theta)$, which is an area nonincreasing map, cannot map surjectively onto the sphere.  

This leads us to make the following conjecture:

\begin{conj} \label{Conj-MinA}
$MinA(\Sph^2\times_{f_j}\Sph^1)=4\pi$ for Example~\ref{ex-sequence}.   
\end{conj}

Recall that
the strong maximum principal states that if there is a smooth set $K$ and an open domain $U$ about $p$ such that
$\Sigma \cap U \subset K$ and $p\in \Sigma \cap \partial K$ then $H(\partial K, p, \nu)= 0$ implies $\Sigma \cap U \subset \partial K$.   
This gives us the following remarks which might be useful to those trying to prove  Conjecture~\ref{Conj-MinA}.

\begin{rmrk} \label{rmrk-phi-min}
We can apply this strong maximum principal to study surfaces which touch $\partial K=\Sigma_{\varphi=\varphi_0}$ from one side to conclude that $\varphi$ cannot achieve a strong local max or min on any $g_j$-minimal surface.    Since we know surfaces of constant $\varphi$,  $\Sigma_{\varphi=\varphi_0}$, have area $\ge 4\pi$, any $g_j$-minimal surface of area $<4\pi$ must wrap around the ${\mathbb{S}}^1$ fibre taking all values of $\varphi$. 
\end{rmrk}

\begin{rmrk} \label{rmrk-theta-min}
We can apply this strong maximum principal to study surfaces which touch $\partial K=\Sigma_{\theta=\theta_0}$ from one side to conclude that $\theta$ cannot achieve a strong local max or min away from the poles where $r=0,\pi$ on any $g_j$-minimal surface.  
Since we know surfaces of constant $\theta$,  $\Sigma_{\theta=\theta_0}$, have area $\ge 4\pi$, any $g_j$-minimal surface of area $<4\pi$ must wrap around the sphere taking on all values of $\theta$.   Be sure to prove this carefully if you wish to use this remark.
\end{rmrk} 

\bibliographystyle{plain}
\bibliography{2023.bib}

\end{document}